\documentclass[11pt]{amsart}
\usepackage{amssymb}

\newcommand{\Sym}{{\mathfrak{S}}}
\newcommand{\Z}{{\mathbb{Z}}}

\newcommand{\N}{{\mathbb{N}}}
\newcommand{\R}{{\mathbb{R}}}

\newcommand{\bT}{{\mathbb{T}}}
\newcommand{\ba}{{\mathbf{a}}}
\newcommand{\bD}{{\mathbf{\Delta}}}
\newcommand{\cO}{{\mathcal{O}}}
\newcommand{\cD}{{\mathcal{D}}}
\newcommand{\bH}{{\mathcal{H}}}
\newcommand{\cC}{{\mathfrak{T}}}
\newcommand{\cL}{{\mathcal{L}}}
\newcommand{\cR}{{\mathcal{R}}}
\newcommand{\cLR}{{\mathcal{LR}}}
\newcommand{\fC}{{\mathfrak{C}}}
\newcommand{\fI}{{\mathfrak{I}}}
\newcommand{\fX}{{\mathfrak{X}}}
\newcommand{\fY}{{\mathfrak{Y}}}
\newcommand{\fS}{{\mathfrak{S}}}
\newcommand{\fR}{{\mathfrak{R}}}
\newcommand{\fp}{{\mathfrak{p}}}
\newcommand{\fs}{{\mathfrak{s}}}
\newcommand{\ft}{{\mathfrak{t}}}

\newcommand{\Irr}{{\operatorname{Irr}}}
\renewcommand{\leq}{\leqslant}
\renewcommand{\geq}{\geqslant}
\renewcommand{\atop}[2]{\genfrac{}{}{0pt}{}{#1}{#2}}

\newtheorem{thm}{Theorem}[section]
\newtheorem{lem}[thm]{Lemma}
\newtheorem{cor}[thm]{Corollary}
\newtheorem{prop}[thm]{Proposition}
\newtheorem{conj}[thm]{Conjecture}

\theoremstyle{definition}
\newtheorem{exmp}[thm]{Example}
\newtheorem{defn}[thm]{Definition}

\theoremstyle{remark}
\newtheorem{rem}[thm]{Remark}

\begin{document}

\title{Lusztig's $a$-function in type $B_n$ in the asymptotic case}
\author{Meinolf Geck and Lacrimioara Iancu}
\address{M.G.: Institut Camille Jordan, Universit\'e Lyon 1, 
21 av Claude Bernard, 69622 Villeurbanne cedex, France}

\curraddr{Department of Mathematical Sciences, King's College,
Aberdeen University, Aberdeen AB24 3UE, Scotland, U.K.}

\email{m.geck@maths.abdn.ac.uk}

\address{L.I.: IGAT, EPFL, B\^atiment de chimie, 1015 Lausanne, Switzerland}

\curraddr{Department of Mathematical Sciences, King's College,
Aberdeen University, Aberdeen AB24 3UE, Scotland, U.K.}

\email{l.iancu@maths.abdn.ac.uk}

\date{April, 2005}
\subjclass[2000]{Primary 20C08; Secondary 20G40}

\dedicatory{To George Lusztig on his $60$th birthday}

\begin{abstract} 
In this paper, we study Lusztig's $\ba$-function for a Coxeter group 
with unequal parameters. We determine that function explicitly in the 
``asymptotic case'' in type $B_n$, where the left cells have been 
determined in terms of a generalized Robinson--Schensted correspondence 
by Bonnaf\'e and the second author. As a consequence, we can show 
that all of Lusztig's conjectural properties (P1)--(P15) hold in this 
case, except possibly (P9), (P10) and (P15). Our methods rely on the 
``leading matrix coefficients'' introduced by the first author. We also
interprete the ideal structure defined by the two-sided cells in the 
associated Iwahori--Hecke algebra $\bH_n$ in terms of the 
Dipper--James--Murphy basis of $\bH_n$.
\end{abstract}

\maketitle

\pagestyle{myheadings}
                                                                               
\markboth{Geck and Iancu}{Lusztig's $a$-function}

\maketitle

\section{Introduction} \label{sec:intro}

Let $(W,S)$ be a Coxeter system where $W$ is finite.  We shall be interested 
in the Kazhdan--Lusztig cells and Lusztig's $\ba$-function on $W$, which
play an important role in the representation theory of finite reductive
groups. The notions of cells and $\ba$-functions are defined in terms of 
the Iwahori--Hecke algebra associated with $W$. Originally, Kazhdan--Lusztig
\cite{KL} and Lusztig \cite{Lu1} only considered the case of a one-parameter 
Iwahori--Hecke algebra; subsequently, the theory has been extended to the 
case of unequal parameters by Lusztig \cite{Lusztig83}, \cite{Lusztig03}. 
However, many results that are known to hold in the equal parameter case
(thanks to a geometric interpretation of the Kazhdan--Lusztig basis) are 
only conjectural in the general case of unequal parameters. A precise
set of conjectures has been formulated by Lusztig in 
\cite[Chap.~14]{Lusztig03}, (P1)--(P15).  (We recall these conjectures
in Section~2.)

The aim of this paper is to determine Lusztig's $\ba$-function explicitly
in the case where $W=W_n$ is of type $B_n$ and the parameters satisfy
the ``asymptotic'' conditions in Bonnaf\'e--Iancu \cite{BI}. As an 
application, we show that all of the conjectures in 
\cite[Chap.~14]{Lusztig03} hold in this case, except possibly (P9), (P10) 
and (P15)\footnote{In a subsequent paper \cite{my05}, using completely 
different methods, the first author shows that (P9), (P10) and a weak 
version of (P15) also hold.  Thus, eventually, (P1)--(P14) and a weak 
version of (P15) are known to hold in the ``asymptotic case'' in type 
$B_n$.}. We also determine the structure of the associated ring $J$. 
Our methods rely on the ``leading matrix coefficients'' introduced by the 
first named author \cite{my02}. It is our hope that similar methods may 
also be applied to other choices of parameters in type $B_n$ where the 
left cell representations are expected to be irreducible. 

In a different direction, we show that the ideal structure defined by the 
two-sided cells in the ``asymptotic'' case in type $B_n$ corresponds 
precisely to the ideal structure given in terms of the Dipper--James--Murphy 
basis \cite{DJM}.

To state our main results more precisely, we have to introduce some notation.
In \cite{Lusztig03}, an Iwahori--Hecke algebra with possibly unequal
parameters is defined with respect to an  integer-valued weight function on
$W$. Following a suggestion of Bonnaf\'e \cite{BI2}, we can slightly modify 
Lusztig's definition so as to include the more general setting in 
\cite{Lusztig83} as well. Let $\Gamma$ be an abelian group (written 
additively) and assume that there is a total order $\leq$ on $\Gamma$ 
compatible with the group structure. (In the setting of \cite{Lusztig03}, 
$\Gamma=\Z$ with the natural order.) 

Let $A={\Z}[\Gamma]$ be the free abelian group with basis $\{e^\gamma
\mid \gamma \in \Gamma\}$. There is a well-defined ring structure on $A$
such that $e^\gamma e^{\gamma'}=e^{\gamma+ \gamma'}$ for all $\gamma,\gamma'
\in \Gamma$.  (Hence, if $\Gamma=\Z$, then $A$ is nothing but the ring of
Laurent polynomials in an indeterminate~$e$.) We write $1=e^0 \in A$.
Given $a\in A$ we denote by $a_\gamma$ the coefficient of $e^\gamma$, so 
that $a= \sum_{\gamma \in \Gamma} a_\gamma e^\gamma$. We denote by 
$A_{\geq 0}$ the set of $\Z$-linear combinations of elements $e^{\gamma}$ 
where  $\gamma\geq 0$. Similarly, we define $A_{>0}$, $A_{\leq 0}$ and 
$A_{<0}$. We say that a function 
\[ L \colon W\rightarrow\Gamma\] 
is a weight function if $L(ww')= L(w)+L(w')$ whenever we have $\ell(ww')=
\ell(w)+\ell(w')$ where $\ell\colon W\rightarrow {\N}$ is the usual length 
function. (We denote $\N=\{0,1,2,\ldots\}$.) We assume throughout that 
$L(s)>0$ for all $s\in S$. Let $\bH=\bH(W,S,L)$ be the generic Iwahori--Hecke 
algebra over $A$ with parameters $\{v_s \mid s\in S\}$ where $v_s:=e^{L(s)}$ 
for $s\in S$.  The algebra $\bH$ is free over $A$ with basis $\{T_w\mid w
\in W\}$, and the multiplication is given by the rule
\[ T_sT_w=\left\{\begin{array}{cl} T_{sw} & \quad \mbox{if $\ell(sw)>\ell
(w)$},\\
T_{sw}+(v_s-v_s^{-1})T_w & \quad \mbox{if $\ell(sw)<\ell(w)$},\end{array}
\right.\]
where $s\in S$ and $w\in W$.
Having fixed a total order on $\Gamma$, we have a corresponding 
Kazhdan--Lusztig basis $\{C_w' \mid w\in W\}$ of $\bH$; we have
\[ C_w'=T_w+\sum_{\atop{y \in W}{y < w}} P_{y,w}^{\,*} \,T_y \in \bH, \]
where $<$ denotes the Bruhat--Chevalley order on $W$ and $P_{y,w}^*\in
A_{<0}$ for all $y<w$ in $W$; see \cite[\S 6]{Lusztig83}. (In the framework 
of \cite{Lusztig03}, the polynomials $P_{y,w}^*$ are denoted $p_{y,w}$ and 
the basis elements $C_w'$ are denoted $c_w$.) Given $x,y\in W$, we write 
\[ C_x' \,C_y'=\sum_{z\in W} h_{x,y,z}\, C_z' \qquad \mbox{where 
$h_{x,y,z} \in A$}.\]
For a fixed $z\in W$, we set 
\[ \ba(z):= \min \{\gamma \geq 0\mid e^\gamma \,h_{x,y,z} \in 
A_{\geq 0} \mbox{ for all $x,y\in W$}\};\]
this is Lusztig's function $\ba \colon W \rightarrow \Gamma$.
(If $\Gamma=\Z$ with its natural order, then this reduces to the 
function defined by Lusztig \cite{Lu1}.)
In Section~2, we recall Lusztig's conjectures concerning the $\ba$-function
and its relation with the pre-order relations $\leq_{\cL}$, $\leq_{\cR}$
and $\leq_{\cLR}$. In the case where $W$ is a Weyl group and $L$ is
constant on $S$, these conjectures are known to hold, thanks to a geometric 
interpretation which yields certain ``positivity properties''; see Lusztig
\cite{Lu1}. In the general case of unequal parameters, it is known 
that these  positivity properties are no longer satisfied.

In this paper, we will be dealing with a Coxeter group of type $B_n$ 
where the parameters are specified as follows.

\begin{exmp} \label{asym} Let $\Gamma$ be any totally ordered abelian
group. Let $W_n$ be a Coxeter group of type $B_n$ ($n \geq 2$), with 
generators,  relations and weight function $L \colon W_n \rightarrow 
\Gamma$ given by the following diagram:
\begin{center}
\begin{picture}(250,50)
\put(  3, 25){$B_n$}
\put(  4, 05){$L$ :}
\put( 40, 25){\circle{10}}
\put( 44, 22){\line(1,0){33}}
\put( 44, 28){\line(1,0){33}}
\put( 81, 25){\circle{10}}
\put( 86, 25){\line(1,0){29}}
\put(120, 25){\circle{10}}
\put(125, 25){\line(1,0){20}}
\put(155, 22){$\cdot$}
\put(165, 22){$\cdot$}
\put(175, 22){$\cdot$}
\put(185, 25){\line(1,0){20}}
\put(210, 25){\circle{10}}
\put( 37, 37){$t$}
\put( 36, 05){$b$}
\put( 76, 37){$s_1$}
\put( 78, 05){$a$}
\put(116, 37){$s_2$}
\put(118, 05){$a$}
\put(203, 37){$s_{n-1}$}
\put(208, 05){$a$}
\end{picture}
\end{center}
where $a,b \in \Gamma$ are such that  
\begin{center} \fbox{$b > (n-1)a>0.$}\end{center}
(Here, $(n-1)a$ means $a+\cdots +a$ in $\Gamma$, with $n-1$ summands.) We 
refer to this hypothesis as the {\bf ``asymptotic case''} in type $B_n$.
Let $\bH_n$ be the corresponding Iwahori--Hecke algebra over 
$A={\Z}[\Gamma]$, where we set
\[ V:=v_t=e^b \qquad \mbox{and}\qquad v:=v_{s_1}=\cdots =v_{s_{n-1}}=e^a.\]
We have the following special case worth mentioning: Let $\Gamma_0=\Z^2$. 
Let $\leq$ be the usual lexicographic order so that $(i,j)<(i',j')$ if 
$i<i'$ or if $i=i'$ and $j<j'$. Then $A_0= {\Z}[\Gamma_0]$ is nothing but 
the ring of Laurent polynomials in two independent indeterminates 
$V_0=e^{(1,0)}$ and~$v_0=e^{(0,1)}$. This is the set-up originally 
considered by Bonnaf\'e--Iancu \cite{BI}; we may refer to this case as 
the {\bf ``generic asymptotic case''} in type $B_n$. 

In Bonnaf\'e--Iancu \cite{BI} (for the ``generic'' case), the left cells 
of $W_n$ are determined explicitly in terms of a generalized 
Robinson--Schensted correspondence. This correspondence associates to each 
element $w\in W_n$ a pair of standard bitableaux of the same shape and total 
size $n$. (By a bitableau, we mean an ordered pair of two tableaux; the
shape of a bitableau is an ordered pair of partitions, that is, a 
bipartition.) Subsequently, Bonnaf\'e \cite{BI2} has shown that these results 
remain valid in the general ``asymptotic case''.  (In Section~5, we recall 
in more detail the main results of \cite{BI}, \cite{BI2}.)
\end{exmp}

Our first main result gives an explicit description of the $\ba$-function.

\begin{thm} \label{main} In the setting of Example~\ref{asym}, let 
$w\in W_n$ and assume that $w$ corresponds to a pair of bitableaux of 
shape $(\lambda_1,\lambda_2)$ by the generalized Robinson--Schensted
correspondance defined in \cite{BI}, where $\lambda_1$ and $\lambda_2$ are 
partitions such that $|\lambda_1|+|\lambda_2|=n$.  Then 
\[ \ba(w)=b\,|\lambda_2|+a\,(n(\lambda_1)+2n(\lambda_2^*)- n(\lambda_2)).\]
Here, $n(\mu)=\sum_{i} (i-1)\mu^{(i)}$ for any partition $\mu=(\mu^{(1)}
\geq \mu^{(2)} \geq \ldots \geq 0)$ and $\mu^*$ denotes the conjugate 
partition. 
\end{thm}

In the ``generic asymptotic case'', the above formula reads
\[ \ba(w)=(|\lambda_2|,\,n(\lambda_1)+2n(\lambda_2^*)- n(\lambda_2)).\]
The proof will be given in Section~5, using the general methods developed
in Section~4. The main ingredients in that proof are Bonnaf\'e's results 
\cite{BI2} on the two-sided cells in $W_n$, and the orthogonal 
representations and leading matrix coefficients introduced in \cite{my02}. 
These are generalizations of the leading coefficients of character values 
considered by Lusztig \cite{Lusztig87}. As an application, we obtain the 
following result. (See Section~5 for the proof.) 

\begin{thm} \label{mainb} In the setting of Example~\ref{asym}, all
the conjectures (P1)--(P15) in \cite[Chap.~14]{Lusztig03} hold except
possibly (P9), (P10) and (P15). The set of ``distinguished involutions''
is given by $\cD=\{z\in W_n \mid z^2=1\}$.
\end{thm}

Thanks to the validity of the properties in Theorem~\ref{mainb}, we can
construct the ring $J$ as explained in \cite[Chap.~18]{Lusztig03}. As an
abelian group, $J$ is free with a basis $\{t_w\mid w\in W_n\}$. The
multiplication is given by 
\[ t_x \cdot t_y=\sum_{z\in W} \gamma_{x,y,z^{-1}} \, t_z\qquad
\mbox{for all $x,y \in W_n$},\]
where $\gamma_{x,y,z^{-1}}\in \Z$ is the constant term of $e^{\ba(z)}\, 
h_{x,y,z} \in A_{\geq 0}$.

\begin{thm} \label{Jring} In the setting of Example~\ref{asym}, we have
an isomorphism of rings $J \cong \bigoplus_{\lambda} M_{d_{\lambda}}(\Z)$,
where $\lambda=(\lambda_1,\lambda_2)$ runs over all bipartitions of $n$ and 
$d_{\lambda}$ is the number of standard bitableaux of shape $\lambda$. 
We have 
\[ \gamma_{x,y,z}=\left\{\begin{array}{cl} \pm 1 & \qquad \mbox{if 
$x\sim_{\cL} y^{-1}$, $y\sim_{\cL} z^{-1}$, $z\sim_{\cL} x^{-1}$},\\
0 & \qquad \mbox{otherwise}.\end{array}\right.\]
\end{thm} 

The proof in Proposition~\ref{strucJ} actually yields an explicit 
isomorphism which shows that $\pm t_w$ ($w\in W_n$) corresponds to a matrix 
unit in $M_{d_{\lambda}}(\Z)$ for some bipartition $\lambda$. Furthermore, 
the signs are interpreted in terms of leading matrix coefficients.

Finally, we show that the Kazhdan--Lusztig basis in the asymptotic case is 
compatible with the basis constructed by Dipper--James--Murphy \cite{DJM}. 
That basis is denoted $\{x_{\fs\ft}\}$ where $(\fs,\ft)$ runs
over all pairs of standard bitableaux of total size $n$ and of the same 
shape; see \cite[Theorem~4.14]{DJM}. Note that the construction of the 
elements $x_{\fs\ft}$ does not rely on the choice of any total order on 
$\Gamma$.  Given a bipartition $\lambda$ of $n$, let $N^\lambda \subseteq 
\bH_n$ be the $A$-submodule spanned by all $x_{\fs\ft}$ where the shape of 
$\fs$ and $\ft$ is a bipartition $\mu$ of $n$ such that $\lambda
\trianglelefteq \mu$. By \cite[Cor.~4.13]{DJM}, $N^\lambda$ is a two-sided 
ideal of $\bH_n$. Now we can state (see the end of Section~5 for the proof):

\begin{thm} \label{murphy} In the setting of Example~\ref{asym}, let 
$\lambda=(\lambda_1,\lambda_2)$ be a bipartition of $n$. Then $N^\lambda$ 
is spanned  by the basis elements $C_w'$ where $w\in W_n$ corresponds, 
via the generalized Robinson--Schensted correspondence, to a bitableau of 
shape $\nu=(\nu_1, \nu_2)$ such that $(\lambda_1,\lambda_2) 
\trianglelefteq (\nu_2,\nu_1^*)$.
\end{thm}

This paper is organized as follows. In Section~2, we recall the basic
definitions concerning the Kazhdan--Lusztig pre-orders on a finite Coxeter 
group $W$ and state Lusztig's conjectures (P1)--(P15), following
\cite[Chap.~14]{Lusztig03}. 

In Section~3, we deal with the leading coefficients of the matrix 
representations of the Iwahori--Hecke algebra associated with $W$ and 
show that, under suitable hypotheses, these leading coefficients can 
be used to detect left, right and two-sided cells. This is an elaboration, 
with some refinements, of the ideas in \cite{my02}. 

In Section~4, we present some criteria and tools for attacking Lusztig's
conjectures. It is our hope that these methods will also be applicable to 
other situations where the left cell representations are expected to be 
irreducible. 

In Section~5, we show that the hypotheses required for the criteria 
in Section~4 are all satisfied for the ``asymptotic case'' in type $B_n$. 
This heavily relies on the fact that Hoefsmit's \cite{Hoefs} matrix 
representations in type $B_n$ are ``orthogonal representations'' in the 
sense of \cite{my02}; hence the theory of leading coefficients and the 
results in Sections~3 and~4 can be applied in this case.

\section{Left cells and Lusztig's conjectures} \label{a-func}

We keep the basic set-up introduced in Section~1 where $W$ is a finite 
Coxeter group and $\bH$ is the corresponding Iwahori--Hecke algebra over 
$A$, with parameters $\{v_s \mid s\in S\}$ where $v_s=e^{L(s)}$ and
$L(s)>0$ for all $s\in S$.

Since we will be dealing 
with $\ba$-invariants of elements in $W$ and of irreducible representations,
it will be technically more convenient to work with a slightly different 
version of the Kazhdan--Lusztig basis of $\bH$. (The reasons can be seen, 
for example, in \cite[Chap.~18]{Lusztig03}.) For any $a\in A$, we define 
$\bar{a}:=\sum_{\gamma \in \Gamma} a_\gamma e^{-\gamma}$. Then we have 
a unique ring involution $j \colon \bH \rightarrow \bH$ such that 
$j(a)=\bar{a}$ for $a\in A$ and $j(T_w)=\varepsilon_wT_{w}$ for $w\in W$, 
where we set $\varepsilon_w=(-1)^{\ell(w)}$. As in \cite[\S 6]{Lusztig83}, 
we set $C_w:=\varepsilon_w j(C_w')$. Then we have 
\[ C_w=T_w+\sum_{\atop{y\in W}{y<w}} \varepsilon_y\,\varepsilon_w\,
\,\overline{P}_{y,w}^{\,*} \,T_y  \qquad \mbox{for all $w\in W$}.\]
The multiplication rule now reads:
\[ C_x\, C_y=\sum_{z\in W} \varepsilon_x\, \varepsilon_y \, \varepsilon_z
\, h_{x,y,z}\,C_z \qquad \mbox{for any $x,y \in W$}.\]
For example, if $x=s\in S$, we have (see \cite[\S 6]{Lusztig83}): 
\[ C_s\,C_y = \left\{\begin{array}{ll} \displaystyle{C_{sy}-
\sum_{\atop{z \in W}{sz<z<y}} \varepsilon_y\, \varepsilon_z\,M_{z,y}^s \,
C_z} &\quad \mbox{if $sy>y$},\\-(v_s+v_s^{-1})\,C_y &\quad \mbox{if $sy<y$},
\end{array}\right.\]
where $M_{z,y}^s \in A$ is determined as in \cite[\S 3]{Lusztig83}. Note 
that we have $\overline{M}_{z,y}^s=M_{z,y}^s$. 

Throughout this paper, we will make use of another important feature of
Iwahori--Hecke algebras, namely, the fact that these algebras carry
a natural {\em symmetrizing trace}. Indeed, consider the linear map 
$\tau \colon \bH \rightarrow A$ defined by $\tau(T_1)=1$ and 
$\tau(T_w)=0$ for $1\neq w\in W$. Then we have 
\[ \tau(T_wT_{w'})=\left\{\begin{array}{cl} 1 & \qquad 
\mbox{if $w'=w^{-1}$},\\ 0 & \qquad \mbox{if $w' \neq w^{-1}$};
\end{array}\right.\]
see \cite[\S 8.1]{ourbuch}. Thus, $\tau$ is a {\em symmetrizing trace} 
on $\bH$. In the following discussion, we shall also need the basis of $\bH$ 
which is dual to the basis $\{C_w\}$ with respect to the symmetrizing 
trace $\tau$. For any $y\in W$ we set 
\[ D_y:=T_y+\sum_{\atop{w \in W}{y <w}} \overline{P}_{ww_0,yw_0}^{\,*}
\,T_w\in \bH,\]
where $w_0\in W$ is the unique element of maximal length in $W$. Then we have
\[ \tau(C_wD_{y^{-1}})=\left\{\begin{array}{cl} 1 & \quad \mbox{if $y=w$},\\
0 & \quad \mbox{if $y \neq w$}.\end{array}\right.\]
(See \cite[2.4]{my02}; see also \cite[Prop.~11.5]{Lusztig03} where the 
analogous statement is proved for the $C'$-basis.)
This immediately yields the following result:

\begin{cor} \label{cor01a} For any $z\in W$, we have 
\[ \ba(z)=\min\{\gamma \geq 0\mid e^\gamma \, \tau(C_xC_y
D_{z^{-1}})\in A_{\geq 0} \mbox{ for all $x,y\in W$}\}.\]
\end{cor}

(Indeed, just note that $\tau(C_xC_yD_{z^{-1}})=\varepsilon_x\,\varepsilon_y
\,\varepsilon_z\,h_{x,y,z}$.)

\medskip
We recall the definition of the left cells of $W$ and the corresponding 
left cell representations of $\bH$ (see \cite{Lusztig83} or \cite{Lusztig03}).

We write $z \leftarrow_{\cL} y$ if there exists some $s\in S$ such that 
$h_{s,y,z} \neq 0$, that is, $C_z'$ occurs in $C_s'\, C_y'$ (when 
expressed in the $C'$-basis) or, equivalently, $C_z$ occurs in $C_s\, C_y$ 
(when expressed in the $C$-basis). Let $\leq_{\cL}$ be the pre-order 
relation on $W$ generated by $\leftarrow_{\cL}$, that is, we have $z 
\leq_{\cL} y$ if there exist elements $z=z_0, z_1,\ldots,z_k=y$ such that 
$z_{i-1} \leftarrow_{\cL} z_i$ for $1\leq i\leq k$.  The equivalence
relation associated with $\leq_{\cL}$ will be denoted by $\sim_{\cL}$ and 
the corresponding equivalence classes are called the {\em left cells} of $W$. 

Similarly, we can define a pre-order $\leq_{\cR}$ by considering
multiplication by $C_s'$ on the right in the defining relation. The 
equivalence relation associated with $\leq_{\cR}$ will be denoted by 
$\sim_{\cR}$ and the corresponding equivalence classes are called the 
{\em right cells} of $W$.  We have
\[ x \leq_{\cR} y \quad \Leftrightarrow \quad x^{-1} \leq_{\cL} y^{-1}.\]
This follows by using the anti-automorphism $\flat\colon \bH\rightarrow
\bH$ given by $T_w^\flat=T_{w^{-1}}$; we have $C_w'^\flat=C_{w^{-1}}'$ and
$C_w^\flat=C_{w^{-1}}$ for all $w\in W$; see \cite[5.6]{Lusztig03}. Thus, 
any statement concerning the left pre-order relation $\leq_{\cL}$ has an 
equivalent version for the right pre-order relation $\leq_{\cR}$, via 
$\flat$. 

Finally, we define a pre-order $\leq_{\cLR}$ by the condition 
that $x\leq_{\cLR} y$ if there exists a sequence $x=x_0,x_1,\ldots, x_k=y$ 
such that, for each $i \in \{1,\ldots,k\}$, we have $x_{i-1} \leq_{\cL} x_i$ 
or $x_{i-1}\leq_{\cR}x_i$. The equivalence relation associated with 
$\leq_{\cLR}$ will be denoted by $\sim_{\cLR}$ and the corresponding 
equivalence classes are called the {\em two-sided cells} of $W$.

Each left cell $\fC$ gives rise to a representation of~$\bH$. This is 
constructed as follows (see \cite[\S 7]{Lusztig83}). Let $[\fC]_A$ be an 
$A$-module with a free $A$-basis $\{c_w \mid w \in \fC\}$. Then the action 
of $C_w$ ($w \in W$) on $[\fC]_A$ is given by the above multiplication 
formulas, i.e., we have 
\[ C_w.c_x = \sum_{y\in \fC} \varepsilon_w\, \varepsilon_x\, \varepsilon_y\,
h_{w,x,y} \, c_y \qquad \mbox{for all $x\in \fC$ and $w\in W$}.\]

\begin{rem} \label{otherrep} It is also possible to define a left
cell module using the $C'$-basis. Recall that $C_w=\varepsilon_w j(C_w')$
for all $w\in W$. Let $\fC$ be a left cell of $W$ and let $[\fC]_A'$ be a 
free $A$-module with a basis $\{c_x'\mid x\in \fC\}$. Then we have an
$\bH$-module structure on $[\fC]_A'$ given by the formula
\[ C_w'.c_x'=\sum_{y\in W} h_{w,x,y}\, c_y' \qquad 
\mbox{for all $x\in \fC$ and $w\in W$}.\]
The passage between the two definitions can be performed using the 
$A$-algebra automorphism $\delta \colon \bH\rightarrow \bH$ given by 
$T_s\mapsto -T_s^{-1}$ for $s\in S$. Note that, by the definition of the 
Kazhdan--Lusztig basis, the elements $C_w$ and $C_w'$ are fixed under 
the composition $\delta \circ j$. Hence, we have $C_w'=\varepsilon_w j(C_w)
= \varepsilon_w\delta(C_w)$ for all $w\in W$. This shows that we have an 
isomorphism of $\bH$-modules
\[ {^\delta}[\fC]_A \cong [\fC]_A' \]
where ${^\delta}[\fC]_A$ is the $\bH$-module obtained from $[\fC]_A$
by composing the original action of $\bH$ with $\delta$.  This remark
will play a role in Example~\ref{bn1} below.
\end{rem}

For the convenience of the reader, we restate here Lusztig's conjectures
(P1)--(P15) in \cite[Chap.~14]{Lusztig03} in the general framework involving
a totally ordered abelian group $\Gamma$. For $z\in W$, we define an 
element $\bD(z)\in \Gamma$ and an integer $0\neq n_z\in \Z$ by the condition
\[e^{\bD(z)}P_{1,z}^* \equiv n_z \quad \bmod A_{<0}; \qquad\mbox{see 
\cite[14.1]{Lusztig03}}.\]
Note that $\bD(z)\geq 0$. Furthermore, given $x,y,z\in W$, we define 
$\gamma_{x,y,z^{-1}}\in \Z$ by
\[ \gamma_{x,y,z^{-1}}=\mbox{constant term of $e^{\ba(z)}\, h_{x,y,z} 
\in A_{\geq 0}$}.\]

\begin{conj}[Lusztig \protect{\cite[14.2]{Lusztig03}}] \label{Pconj}
Let $\cD=\{z\in W \mid \ba(z)=\bD(z)\}$.  Then the  following properties
hold.
\begin{itemize}
\item[\bf P1.] For any $z\in W$ we have $\ba(z)\leq \bD(z)$.
\item[\bf P2.] If $d \in \cD$ and $x,y\in W$ satisfy $\gamma_{x,y,d}\neq 0$,
then $x=y^{-1}$.
\item[\bf P3.] If $y\in W$, there exists a unique $d\in \cD$ such that
$\gamma_{y^{-1},y,d}\neq 0$.
\item[\bf P4.] If $z'\leq_{\cL\cR} z$ then $\ba(z')\geq \ba(z)$. Hence, if
$z'\sim_{\cL\cR} z$, then $\ba(z)=\ba(z')$.
\item[\bf P5.] If $d\in \cD$, $y\in W$, $\gamma_{y^{-1},y,d}\neq 0$, then
$\gamma_{y^{-1},y,d}=n_d=\pm 1$.
\item[\bf P6.] If $d\in \cD$, then $d^2=1$.
\item[\bf P7.] For any $x,y,z\in W$, we have $\gamma_{x,y,z}=\gamma_{y,z,x}$.
\item[\bf P8.] Let $x,y,z\in W$ be such that $\gamma_{x,y,z}\neq 0$. Then
$x\sim_{\cL} y^{-1}$, $y \sim_{\cL} z^{-1}$, $z\sim_{\cL} x^{-1}$.
\item[\bf P9.] If $z'\leq_{\cL} z$ and $\ba(z')=\ba(z)$, then $z'\sim_{\cL}z$.
\item[\bf P10.] If $z'\leq_{\cR} z$ and $\ba(z')=\ba(z)$, then $z'\sim_{\cR}z$.
\item[\bf P11.] If $z'\leq_{\cL\cR} z$ and $\ba(z')=\ba(z)$, then
$z'\sim_{\cL\cR}z$.
\item[\bf P12.] Let $I\subset S$ and $W_I$ be the parabolic subgroup
generated by $I$. If $y\in W_I$, then $\ba(y)$ computed in terms of $W_I$
is equal to $\ba(y)$ computed in terms of~$W$.
\item[\bf P13.] Any left cell $\fC$ of $W$ contains a unique element
$d\in \cD$. We have $\gamma_{x^{-1},x,d}\neq 0$ for all $x\in \fC$.
\item[\bf P14.] For any $z\in W$, we have $z \sim_{\cL\cR} z^{-1}$.
\item[\bf P15.] If $x,x',y,w\in W$ are such that $\ba(w)=\ba(y)$, then
\[\sum_{y' \in W} h_{w,x',y'}\otimes_{\Z}h_{x,y',y}=
\sum_{y'\in W} h_{x,w,y'}\otimes_{\Z} h_{y',x',y} \quad \mbox{in 
$A \otimes_{\Z} A$}. \]
\end{itemize}
\end{conj}
(The above formulation  of (P15) is taken from Bonnaf\'e \cite{BI2}.)

\begin{rem} \label{invers} For all $x,y,z\in W$, we have
\[ h_{x,y,z}=h_{y^{-1},x^{-1},z^{-1}} \qquad \mbox{and} \qquad
P_{x,y}^*=P_{x^{-1},y^{-1}}^*.\]
Hence, we have $\ba(z)=\ba(z^{-1})$, $n_z=n_{z^{-1}}$, $\bD(z)=
\bD(z^{-1})$, $\cD=\cD^{-1}$.
\end{rem}

\begin{proof} We have already remarked above that there is an 
anti-automorphism $\flat \colon \bH \rightarrow \bH$ such that 
$T_w^\flat=T_{w^{-1}}$ for all $w \in W$. By the argument in 
\cite[5.6]{Lusztig03}, we have $C_w^\flat=C_{w^{-1}}$. This yields all
the above statements.
\end{proof}

If $W$ is the symmetric group $\fS_n$, the above conjectures are all known
to hold; see \cite[Chap.~15]{Lusztig03}\footnote{In a recent preprint
\cite{my05a}, the first named author has given elementary proofs of 
(P1)--(P15) for $W=\fS_n$.}. Hence 
the information about left, right and two-sided cells, as well as the 
$\ba$-function, is rather complete in this case.  We close this section by
summarizing some known results on $\fS_n$.  (This information will be 
needed in the proof of Theorem~1.2; see Section~5.) 

\begin{exmp} \label{symgroup} Let $\fS_n=\langle s_1, \ldots,s_{n-1}\rangle$
be the symmetric group, where $s_i=(i,i+1)$ for $1\leq i \leq n-1$.  The 
diagram is given as follows.
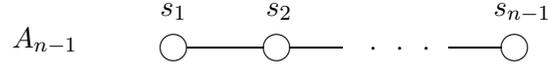
\begin{center}
\begin{picture}(300,30)
\put( 40, 5){$A_{n-1}$}
\put(101, 5){\circle{10}}
\put(106, 5){\line(1,0){29}}
\put(140, 5){\circle{10}}
\put(145, 5){\line(1,0){20}}
\put(175, 2){$\cdot$}
\put(185, 2){$\cdot$}
\put(195, 2){$\cdot$}
\put(205, 5){\line(1,0){20}}
\put(230, 5){\circle{10}}
\put( 96, 17){$s_1$}
\put(136, 17){$s_2$}
\put(222, 17){$s_{n-1}$}
\end{picture}
\end{center}
We consider the abelian group $\Gamma=\Z$ with its natural order and
denote $v:=e^1$. Then $A={\Z}[v,v^{-1}]$ is the ring of Laurent polynomials
in an indeterminate~$v$. Let $L \colon W \rightarrow \Z$ be any 
weight function such that $L(s_i)>0$ for all$~i$. Since all generators are 
conjugate, $L$ takes the same value on each~$s_i$. Thus, we are in the case
of ``equal parameters''.

The classical Robinson--Schensted correspondence associates with
each element $\sigma \in \fS_n$ a pair of standard tableaux $(P(\sigma),
Q(\sigma))$ of the same shape. The tableau $P(\sigma)$ is obtained by
``row-insertion'' of the numbers $\sigma.1,\ldots,\sigma.n$ (in this
order) into an initially empty tableau; the tableau $Q(\sigma)$ ``keeps 
the record'' of the order by which the positions in $P(\sigma)$ have been 
filled; see Fulton \cite[Chap.~4]{Fulton}. For any partition $\nu 
\vdash n$, we set
\[ \fR_\nu:=\{\sigma \in \fS_n \mid \mbox{ $P(\sigma)$, $Q(\sigma)$ have
shape $\nu$}\}.\]
Thus, we have $\fS_n=\coprod_{\nu \vdash n} \fR_\nu$. Then the 
following hold.
\begin{itemize}
\item[(a)] {\em For a fixed standard tableau $T$, the set $\{\sigma \in 
\fS_n \mid Q(\sigma)=T\}$ is a left cell of $\fS_n$ and $\{\sigma \in
\fS_n \mid P(\sigma)=T\}$ is a right cell of $\fS_n$. Furthermore, all 
left cells and all right cells arise in this way}.
\end{itemize}
This was first proved by Kazhdan--Lusztig \cite[\S 4]{KL}; for a more 
direct and elementary proof, see Ariki \cite{Ar}.

\begin{itemize}
\item[(b)] {\em The sets $\fR_\nu ,\; \nu \vdash n$ are precisely the 
two-sided cells of $\fS_n$}. 
\end{itemize}
This is seen as follows. First note that the statements (P1)--(P15)
in Conjecture~\ref{Pconj} are known to hold for $W=\fS_n$ (see
\cite[Chap.~15]{Lusztig03} and the references there). Now (P4), (P9), (P10)
imply that $x,y \in \fS_n$ lie in the same two-sided cell if and only if
there exists a sequence of elements $x=x_0,x_1,\ldots,x_k=y$ in $\fS_n$
such that, for each $i$, we have $x_{i-1} \sim_{\cL} x_i$ or $x_{i-1}
\sim_{\cR} x_i$. Now (b) follows from (a).
\begin{itemize}
\item[(c)] {\em For any $\nu\vdash n$, we have $\sigma_{\nu^*} \in 
\fR_\nu$, where $\sigma_{\nu^*}$ is the longest element in the Young subgroup 
$\fS_{\nu^*} \subseteq \fS_n$ and $\nu^*$ denotes the conjugate partition}.
\end{itemize}
This is a purely combinatorial exercice: it is enough to apply the
Robinson--Schensted correspondence to the element $\sigma_{\nu*}$ and to
verify that the corresponding tableaux have shape $\nu$.
\begin{itemize}
\item[(d)] {\em If $\sigma \in \fR_\nu$, then $\ba(\sigma)=n(\nu)$, 
where $n(\nu)$ is defined as in Theorem~1.2}.
\end{itemize}
This is seen as follows. Again, we use the fact that (P1)--(P15) in
Conjecture~\ref{Pconj} hold for $W=\fS_n$. Since, by (P4), the $\ba$-function
is constant on the two-sided cells, (c) shows that it is enough to 
compute $\ba(\sigma_{\nu})$ for any $\nu\vdash n$. But, since $\sigma_\nu$
is the longest element in a parabolic subgroup, we have $\ba(\sigma_\nu)=
\ell(\sigma_\nu)$ by (P12) and \cite[13.8]{Lusztig03}. It remains to note 
that $\ell(\sigma_\nu)=n(\nu^*)$.
\begin{itemize}
\item[(e)] {\em If $\sigma \in \fR_{\nu}$ and $\sigma' \in \fR_{\nu'}$ are 
such that $\sigma \leq_{\cLR} \sigma'$, then we have $\nu 
\trianglelefteq \nu'$, where $\trianglelefteq$ denotes the dominance
order. This means that 
\[\sum_{i=1}^k \nu^{(i)} \leq \sum_{i=1}^k {\nu'}^{(i)} \quad
\mbox{for all $k>0 $},\]
where $\nu^{(i)}$ and ${\nu'}^{(i)}$ are the parts of $\nu$ and $\nu'$, 
respectively}. 
\end{itemize}
This follows from a result of Lusztig--Xi \cite[3.2]{LuXi}; see 
Du--Parshall--Scott \cite[2.13.1]{DPS} and the references there. (Note that,
here again, (P1)--(P15) are used.)
\end{exmp}

\section{Leading matrix coefficients} \label{lead}

We now recall the basic facts concerning the leading matrix coefficients
introduced in \cite{my02}. We extend scalars from $\Z$ to $\R$ and consider
the group algebra ${\R}[\Gamma]$. Since $\Gamma$ is totally ordered, 
${\R}[\Gamma]$ is an integral domain; let $K$ be its field of fractions. 
We define $\fI_{>0}\subset {\R}[\Gamma]$ to be the set of all $f \in 
{\R}[\Gamma]$ such that
\[ f=1+\mbox{$\R$-linear combination of elements of $e^\gamma$ where 
$\gamma>0$}.\]
Note that $\fI_{>0}$ is multiplicatively closed. Furthermore, every
element $x\in K$ can be written in the form
\[ x=r_x\,e^{\gamma_x}f/g\qquad \mbox{where $r_x \in \R$, $\gamma_x \in 
\Gamma$ and $f,g\in \fI_{>0}$};\]
note that, if $x\neq 0$, then $r_x$ and $\gamma_x$ indeed are 
{\em uniquely determined} by $x$; if $x=0$, we have $r_0=0$ and we set 
$\gamma_0:=+\infty$ by convention. We set 
\[{\cO}:=\{x \in K \mid \gamma_x \geq 0\} \qquad \mbox{and}\qquad
{\fp}:=\{x \in K \mid \gamma_x >0\}.\]
Then it is easily verified that $\cO$ is a valuation ring in $K$, that is,
$\cO$ is a subring of $K$ such that, for any $0 \neq x \in K$, we have 
$x \in \cO$ or $x^{-1}\in \cO$. Furthermore, $\cO$ is a local ring with 
maximal ideal $\fp$. The group of units in $\cO$ is given by
\[ {\cO}^\times=\{x \in \cO \mid r_x\neq 0, \gamma_x=0\}.\]
Note that we have 
\begin{align*}
\cO \cap {\R}[\Gamma]&={\R}[\Gamma]_{\geq 0}:=\langle e^\gamma
\mid \gamma\geq 0 \rangle_{\R},\\
\fp\cap {\R}[\Gamma]&={\R}[\Gamma]_{>0}:=\langle e^\gamma\mid
\gamma>0\rangle_{\R}.
\end{align*}
We have a well-defined $\R$-linear ring homomorphism $\cO \rightarrow \R$ 
with kernel $\fp$. The image of $x\in \cO$ in $\R$ is called the 
{\em constant term} of $x$. Thus, the constant term of $x$ is $0$ if 
$x\in \fp$; the constant term equals $r_x$ if $x\in \cO^\times$. 

Extending scalars from $A$ to $K$, we obtain a finite dimensional 
$K$-algebra $\bH_K$, with basis $\{T_w\mid w\in W\}$ and multiplication
as specified in Section~1. We have:

\begin{rem} \label{semis} The algebra $\bH_K$ is split semisimple and 
abstractly isomorphic to the group algebra of $W$ over $K$.
\end{rem}

\begin{proof} Since the situation here is somewhat more general than usual,
let us indicate the main ingredients. To show that $\bH_K$ is semisimple,
we use the $\R$-linear ring homomorphism $\theta \colon {\R}[\Gamma]
\rightarrow \R$ such that $\theta(e^\gamma)=1$ for all $\gamma \in \Gamma$. 
By extension of scalars, we obtain $\R \otimes_{{\R}[\Gamma]} \bH 
\cong {\R}[W]$, the group algebra of $W$ over $\R$. Since the latter 
algebra is known to be semisimple, a standard argument (using Tits' 
Deformation Theorem) shows that $\bH_K$ must be semisimple, too. (See, 
for example, \cite[7.4.6 and 8.1.7]{ourbuch}.) But then it is 
also known that $\bH_K$ is split and abstractly isomorphic to $K[W]$; see 
\cite[9.3.5 and 9.3.9]{ourbuch} and the references there. 
\end{proof}

Let $\Irr(\bH_K)$ be the set of irreducible characters of $\bH_K$. We 
write this set in the form 
\[ \Irr(\bH_K)=\{\chi_\lambda \mid \lambda\in \Lambda\},\]
where $\Lambda$ is some finite indexing set. The algebra $\bH_K$ is 
{\em symmetric} with respect to the trace function $\tau \colon \bH_K
\rightarrow K$ defined by $\tau(T_1)=1$ and $\tau(T_w)=0$ for $1\neq w\in W$
(see Section~2). The fact that $\bH_K$ is split semisimple yields that 
\[ \tau=\sum_{\lambda\in \Lambda} \frac{1}{c_\lambda} \, \chi_\lambda
\qquad \mbox{where $0 \neq c_\lambda \in {\R}[\Gamma]$}.\]
The elements $c_\lambda$ are called the {\em Schur elements}. By 
\cite[8.1.8]{ourbuch}, we have $c_\lambda=P_{W,L}/D_\lambda$ where 
$P_{W,L}= \sum_{w\in W} e^{2L(w)}$ is the Poincar\'e polynomial of $W,L$ 
and $D_\lambda$ is the ``generic degree'' associated with $\chi_\lambda$. 
We can write
\[ c_\lambda=r_\lambda \, e^{-2\alpha_\lambda}\,f_\lambda\qquad 
\mbox{where $r_\lambda \in \R_{>0}$, $f_\lambda \in {\fI}_{>0}$ and 
$\alpha_\lambda \geq 0$}.\] 
The element $\alpha_\lambda\in \Gamma$ is called the {\em generalized 
$a$-invariant} of $\chi_\lambda$; see \cite[\S 3]{my02}. (Note that the 
notation in [{\em loc.\ cit.}] has to be adapted to the present setting 
where we write the elements of ${\R}[\Gamma]$ exponentially.)

By \cite[Prop.~4.3]{my02}, every $\chi_\lambda$ is afforded by a so-called
orthogonal representation.  This means that there exists a matrix 
representation $\fX_\lambda \colon \bH_K\rightarrow M_{d_\lambda}(K)$ with 
character $\chi_\lambda$ and an invertible diagonal matrix $P \in 
M_{d_\lambda}(K)$ such that the following conditions hold:
\begin{itemize}
\item[(O1)] We have $\fX_\lambda(T_{w^{-1}})=P^{-1} \cdot \fX_\lambda
(T_w)^{\text{tr}}\cdot P$ for all $w\in W$, and 
\item[(O2)] the diagonal entries of $P$ lie in $\fI_{>0}$.
\end{itemize}
This has the following consequence. Let $\lambda \in \Lambda$
and $1\leq i,j \leq d_\lambda$. For any $h\in \bH_K$, we denote by 
$\fX_\lambda^{ij}(h)$ the $(i,j)$-entry of the matrix $\fX_\lambda(h)$. 
Then, by \cite[Theorem~4.4 and Remark~4.5]{my02}, we have
\[ e^{\alpha_\lambda} \fX_\lambda^{ij}(T_w) \in \cO,  \qquad 
e^{\alpha_\lambda} \fX_\lambda^{ij}(C_w) \in \cO,  \qquad 
e^{\alpha_\lambda} \fX_\lambda^{ij}(D_w) \in \cO\]
for any $w\in W$ and 
\[ e^{\alpha_\lambda} \fX_\lambda^{ij}(T_w) \equiv e^{\alpha_\lambda} 
\fX_\lambda^{ij}(C_w) \equiv e^{\alpha_\lambda} 
\fX_\lambda^{ij}(D_w) \bmod\fp.\]
Hence, the above three elements of $\cO$ have the same constant term 
which we write as $\varepsilon_w\,c_{w,\lambda}^{ij}$. The constants 
$c_{w,\lambda}^{ij}\in \R$ are called the {\em leading matrix coefficients} 
of $\fX_\lambda$. By \cite[Theorem~4.4]{my02}, 
these coefficients have the following property:
\begin{alignat*}{2}
c_{w,\lambda}^{ij}&=c_{w^{-1},\lambda}^{ji} &&\qquad 
\mbox{for all $w \in W$},\\ c_{w,\lambda}^{ij} &\neq 0 &&\qquad 
\mbox{for some $w\in W$}.
\end{alignat*}
Furthermore, we have 
\begin{align*}
\alpha_\lambda &=\min\{ \gamma \geq 0 \mid e^\gamma\, \fX_\lambda^{ij}(T_w)
\in \cO \mbox{ for all $w \in W$ and $1\leq i,j\leq d_\lambda$}\}\\
&=\min\{ \gamma \geq 0 \mid e^\gamma\, \fX_\lambda^{ij}(C_w)
\in \cO \mbox{ for all $w \in W$ and $1\leq i,j\leq d_\lambda$}\}\\
&=\min\{ \gamma \geq 0 \mid e^\gamma\, \fX_\lambda^{ij}(D_w)
\in \cO \mbox{ for all $w \in W$ and $1\leq i,j\leq d_\lambda$}\}.
\end{align*}
The leading matrix coefficients satisfy the following {\em Schur relations}.
Let $\lambda,\mu \in \Lambda$, $1\leq i,j \leq d_\lambda$ and $1 \leq k,l
\leq d_\mu$; then 
\[ \sum_{w \in W} c_{w,\lambda}^{ij}c_{w,\mu}^{kl}=\left\{\begin{array}{cl}
\delta_{ik}\,\delta_{jl}\, r_\lambda & \mbox{if $\lambda=\mu$},
\\ 0 & \mbox{if $\lambda\neq \mu$}; \end{array}\right.\]
see \cite[Theorem~4.4]{my02}. Since $|W|=\sum_{\lambda \in \Lambda} 
d_\lambda^2$, we can invert the above relations and obtain another set 
of relations (analogous to the ``second'' orthogonality relations for
the characters of a finite group): For any $y,w\in W$ we have
\[ \sum_{\lambda\in \Lambda} \sum_{i,j=1}^{d_\lambda} 
\frac{1}{r_\lambda} \, c_{y,\lambda}^{ij}c_{w,\lambda}^{ij}=
\left\{\begin{array}{cl} 1 & \mbox{if $y=w$}, \\ 0 & \mbox{if $y\neq w$}.
\end{array}\right.\]
The above relations immediately imply that $W=\bigcup_{\lambda 
\in \Lambda} \cC_\lambda$, where we set 
\[ \cC_\lambda:=\{w\in W\mid c_{w,\lambda}^{ij} \neq 0
\mbox{ for some $1\leq i,j \leq d_\lambda$\}}.\]
The leading matrix coefficients are related to the left cells of $W$ by
the following result. Recall that, given a left cell $\fC$, we have a 
corresponding left cell module $[\fC]_A$. Extending scalars from $A$ to
$K$, we obtain an $\bH_K$-module $[\fC]_K$. We denote by $\chi_{\fC}$ the
character of $[\fC]_K$. Then we have:

\begin{prop}[See \protect{\cite[Prop.~4.7]{my02}}] \label{mylem}
Let $\lambda \in \Lambda$ and $\fC$ be a left cell in $W$.  Denote by 
$[\chi_{\fC}:\chi_\lambda]$ the multiplicity of $\chi_\lambda$ in 
$\chi_{\fC}$.  Then we have 
\[ \sum_{k=1}^{d_\lambda} \sum_{y\in \fC} (c_{y,\lambda}^{ik})^2=
[\chi_{\fC}:\chi_\lambda] \, r_\lambda, \qquad \mbox{for any $1\leq
i \leq d_\lambda$}.\]
In particular, if $w\in \cC_\lambda$, then $\chi_\lambda$ occurs with 
non-zero multiplicity in $\chi_{\fC}$, where $\fC$ is the left cell 
containing $w$.
\end{prop}

\begin{proof} The formula is proved in [{\em loc.\ cit.}]. Now fix 
$\lambda \in \Lambda$ and let $w\in \cC_\lambda$; then $c_{w,
\lambda}^{ij}\neq 0$ for some $1\leq i,j\leq d_\lambda$. Let $\fC$ be
the left cell containing $w$. Now all terms on the left hand side of
the formula are non-negative, and the term corresponding to $y=w$
and $k=j$ is strictly positive. Hence the left hand side is non-zero
and so $[\chi_{\fC}:\chi_\lambda] \neq 0$.
\end{proof}

The Schur relations lead to particularly strong results when some additional
hypotheses are satisfied. These are isolated in the following definition.

\begin{defn} \label{def01} We say that $\bH$ is {\em integral} if 
\[ c_{w,\lambda}^{ij} \in \Z \qquad \mbox{for all $w\in W$, $\lambda 
\in \Lambda$ and $1\leq i ,j \leq d_\lambda$}.\]
Furthermore, recall that $c_\lambda= r_\lambda e^{-2\alpha_\lambda}
f_\lambda$, where $r_\lambda$ is a positive real number and $f_\lambda
\in \fI_{>0}$. We say that $\bH$ is {\em normalized} if 
\[ r_\lambda=1 \qquad \mbox{for all $\lambda \in \Lambda$}.\]
\end{defn}

\begin{rem} \label{rem4a} Suppose that $\bH$ is normalized and that 
Lusztig's conjectures (P1)--(P15) in \cite[Chap.~14]{Lusztig03} are 
satisfied. Then we necessarily have that $\chi_{\fC} \in \Irr(\bH_K)$ 
for all left cells $\fC$ in $W$; see \cite[Cor.~4.8]{my04b}. Thus,
the conditions in Definition~\ref{def01} should be considered as rather 
severe restrictions on the structure of $\bH$.
\end{rem}

Since the Schur elements $c_\lambda$ are known in all cases (see the
appendices of Carter \cite{Ca2} and Geck--Pfeiffer \cite{ourbuch}, for 
example) the condition that $\bH$ be {\em normalized} is rather 
straightforward to verify. The condition that $\bH$ be {\em integral} is 
more subtle. Let us describe here a convenient setting in which this 
condition can be dealt with.  Let $\lambda \in \Lambda$ and assume that 
we have a matrix representation $\fY_\lambda \colon \bH_K \rightarrow 
M_{d_\lambda}(K)$ affording $\chi_\lambda$ such that $\fY_\lambda$ is 
{\em integral}, in the sense that 
\[ \fY_\lambda(T_w) \in M_d(A) \qquad \mbox{for all $w\in W$}.\]
Furthermore, assume that there is a {\em symmetric}  matrix 
$\Omega=(\omega_{ij}) \in M_{d_\lambda}(A)$ satisfying the following two 
conditions:
\begin{itemize}
\item[(F1)] We have $\fY_\lambda(T_{w^{-1}})=\Omega^{-1}\cdot 
\fY_\lambda(T_w)^{\text{tr}} \cdot \Omega$ for all $w\in W$, and 
\item[(F2)] all principal minors of $\Omega$ lie in $1+A_{>0}$.
\end{itemize}
Note that $1+A_{>0}=A \cap \fI_{>0}$.

\begin{lem} \label{critnorm} In the above setting, there is an orthogonal 
representation $\fX_\lambda \colon \bH_K \rightarrow M_{d_\lambda}(K)$ 
affording $\chi_\lambda$ such that the corresponding leading coefficients
satisfy $c_{w,\lambda}^{ij} \in \Z$ for all $w\in W$ and $1\leq i,j\leq 
d_\lambda$. We have $\fX_\lambda(T_w)\in M_{d_\lambda}(K_0)$, where $K_0$
is the field of fractions of $A$.
\end{lem}

\begin{proof} Let $R:=\{f/g \mid f \in A, g \in 1+A_{>0}\} \subseteq K$, 
a subring of $K$.

Now we consider the system of equations $\Pi^{\text{tr}} \cdot Z \cdot 
\Pi=\Omega$, where $\Pi=(\pi_{ij})$ is an {\em upper triangular} matrix 
with $1$ on the diagonal (and unknown coefficients above the diagonal) and 
$Z=(z_{ij})$ is a {\em diagonal} matrix (with unknown coefficients on the 
diagonal). It is well-known and easy to see that the above system has a 
unique solution $(\Pi,Z)$, where $z_{ii}\in R^\times$ for all $i$ 
and $\pi_{ij} \in R$ for all $i<j$. (Note that all principal minors of 
$\Omega$ are invertible in $R$.) Since $\Omega$ and $Z$ have the 
same principal minors, there exists some $f \in 1+A_{>0}$ such that 
\[ fz_{ii} \in 1+A_{>0} \subseteq \fI_{>0} \quad \mbox{for all $i$}.\]
We set $P:=fZ$ and define $\fX_\lambda$ by 
\[ \fX_\lambda(T_w):=\Pi\cdot \fY_\lambda(T_w)\cdot \Pi^{-1}
\qquad\mbox{for all $w\in W$}.\]
Then, clearly, $\fX_\lambda$ affords $\chi_\lambda$, and a straightforward 
computation yields 
\[\fX_\lambda(T_{w^{-1}}) = Z^{-1}\cdot \fX_\lambda(T_w)^{\text{tr}} 
\cdot Z= P^{-1}\cdot \fX_\lambda(T_w)^{\text{tr}} \cdot P.\]
So $\fX_\lambda$ and $P$ satisfy the conditions (O1) and (O2). Now
note that $\Pi$ (being triangular with $1$ on the diagonal) is invertible
over $R$; let us denote 
\[ \Pi^{-1}=(\tilde{\pi}_{ij}) \qquad \mbox{where $\tilde{\pi}_{ii}=1$
and $\tilde{\pi}_{ij}\in R$ for all $i,j$}.\]
Hence, for any $1\leq i,j\leq d_\lambda$, we obtain
\[ e^{\alpha_\lambda}\fX_\lambda^{ij}(T_w)=\sum_{k=1}^{d_\lambda} 
\sum_{l=1}^{d_\lambda} {\pi}_{ik}\,\tilde{\pi}_{lj} \Bigl(e^{\alpha_\lambda}
\fY_\lambda^{kl}(T_w) \Bigr)\in R.\]
On the other hand, $\fX_\lambda$ is an orthogonal representation; hence
the above element lies in $\cO \cap R$. Consequently, the corresponding
leading matrix coefficient $c_{w,\lambda}^{ij}$ will lie in $\Z$, as required.
\end{proof}

\begin{exmp} \label{bn} Let $W=W_n$ and $L\colon W_n\rightarrow \Gamma$ 
be as in Example~\ref{asym} (the ``asymptotic case'' in type $B_n$). Let 
$\bH_n$ be the corresponding Iwahori--Hecke algebra and write $\bH_{n,K}=
K \otimes_A \bH_n$.  Then we have a natural parametrization 
\[\Irr(\bH_{n,K})=\{\chi_{\lambda} \mid \lambda \in \Lambda_n\}\]
where $\Lambda_n$ is the set of all bipartitions $\lambda=(\lambda_1,
\lambda_2)$ such that $|\lambda_1|+ |\lambda_2|=n$; see 
\cite[Chap.~5]{ourbuch}. For any $(\lambda_1,\lambda_2) \in \Lambda_n$,
we have
\[r_{(\lambda_1,\lambda_2)}=1 \quad \mbox{and}\quad \alpha_{(\lambda_1,
\lambda_2)}=b\,|\lambda_2|+a\,(n(\lambda_1)+2n(\lambda_2)-
n(\lambda_2^*)),\]
where we set $n(\nu)=\sum_{i} (i-1)\nu^{(i)}$ for any partition $\nu=
(\nu^{(1)} \geq \cdots \geq \nu^{(r)}\geq 0)$ and where $\nu^*$ denotes 
the conjugate partition.  Thus, $\bH_n$ is normalized. (For the proof, 
see \cite[Remark~5.1]{my02} and the references there. Actually, in 
[{\em loc.\ cit.}], we only consider the ``generic asymptotic case''. But 
it is readily checked, using the formula for $c_\lambda$ in 
\cite[p.~447]{Ca2}, that the above formulas also hold in the 
``asymptotic case''. For a weight function with values in $\Z$, explicit 
formulas for $r_{(\lambda_1,\lambda_2)}$ and $\alpha_{(\lambda_1,\lambda_2)}$ 
are given by Lusztig \cite[Prop.~22.14]{Lusztig03}; then it is a purely 
combinatorial exercise to show that, in the case $b>(n-1)a$, Lusztig's 
formulas can be rewritten as above.)

Furthermore, by the discussion in \cite[\S 5]{my02}, $\bH_n$ is integral.  
Since that discussion is somewhat sketchy (and only deals with the 
``generic asymptotic case''), let us add a more rigorous argument, based on 
Dipper--James--Murphy \cite{DJM}, \cite{djm2} and Lemma~\ref{critnorm}. 
That argument also shows that we can work with the field of fractions of
$A$ instead of $K$.

Let $\lambda\in \Lambda_n$ and $\tilde{S}^\lambda$ be the corresponding 
Specht module over $K$, as constructed in \cite[4.19]{DJM}. The modules 
$\tilde{S}^\lambda$ are absolutely irreducible and pairwise non-isomorphic 
\cite[4.22]{DJM}. Let $\bT_\lambda$ be the set of standard bitableaux of
shape $\lambda$. Then $\tilde{S}^\lambda$ has a standard basis $\{e_t\mid
t \in \bT_\lambda\}$ such that the corresponding matrix representation is 
{\em integral} in the above sense. Furthermore, by \cite[\S 5]{DJM}, there 
is a non-degenerate bilinear form $\langle \;,\; \rangle_\lambda$ on 
$\tilde{S}^\lambda$, satisfying the condition $\langle T_sx,x' 
\rangle_\lambda=\langle x,T_{s}x' \rangle_\lambda$ for all $x,x'\in 
\tilde{S}^\lambda$ and all $s\in S$. The Gram matrix of $\langle \;,
\;\rangle_{\lambda}$ with respect to the standard basis of 
$\tilde{S}^\lambda$ has coefficients in $A$ and a non-zero determinant. 
We will have to slightly modify the standard basis of $\tilde{S}^\lambda$ 
in order to make sure that (F1) and (F2) hold for the corresponding matrix 
representation.

Now, by \cite[Theorem~8.11]{DJM}, there is an orthogonal basis $\{f_t
\mid t \in \bT_\lambda\}$ with respect to the above bilinear form; moreover, 
the matrix transforming the standard basis into the orthogonal basis is 
triangular with $1$ on the diagonal\footnote{It is remarked in \cite{DJM} 
that the orthogonal representation with respect to the basis $\{f_t\}$ 
actually coincides with the one defined by Hoefsmit \cite{Hoefs}.}. Using
the recursion formula in \cite[Proposition~3.8]{djm2}, it is straightforward
to show that, for each basis element $f_t$, there exist integers $s_t,
a_{ti}, b_{tj},c_{tk}, d_{tl}\in \Z$ such that 
\[ a_{ti} \geq 0,  \qquad b_{tj} \geq 0, \qquad b+c_{tk}a>0, \qquad
b+d_{tl}a>0, \]
and 
\[ \langle f_t,f_t\rangle_\lambda=e^{2s_ta}\cdot \frac{\prod_i 
(1+e^{2a}+ \cdots +e^{2a_{ti}a})}{\prod_j (1+e^{2a}+\cdots +e^{2b_{tj}a})} 
\cdot \frac{\prod_k \bigl(1+e^{2(b+c_{tk}a)}\bigr)}{\prod_l\bigl(1+
e^{2(b+d_{tl}a)}\bigr)}.\]
The fact that $b+c_{tk}a$ and $b+d_{tl}a$ are strictly positive essentially
relies on the condition that $b>(n-1)a$. Hence, setting 
\[ \tilde{f}_t:=e^{-s_ta} \cdot \Bigl(\prod_j (1+e^{2a}+ \cdots 
+e^{2b_{tj}a})\Bigr)\cdot \Bigl(\prod_l(1+e^{2(b+d_{tl}a)})\Bigr)\cdot f_t,\]
we obtain $\langle \tilde{f}_t,\tilde{f}_t\rangle_\lambda \in 1+A_{>0}$ for 
all $t$. Now let 
\[ \tilde{e}_t=e^{-s_ta}e_t \qquad \mbox{for all $t \in \bT_\lambda$}.\]
Then $\{\tilde{e}_t\mid t \in \bT_\lambda\}$ is a new basis of 
$\tilde{S}^\lambda$. Let $\fY_\lambda$ be the corresponding matrix 
representation and $\Omega$ be the corresponding Gram matrix of $\langle 
\;, \; \rangle_{\lambda}$. Then $\fY_\lambda$ is still integral, $\Omega$ 
has coefficients in $A$ and (F1) holds. Now, the matrix transforming 
$\tilde{e}_t$ to the basis $\tilde{f}_t$ is triangular with quotients of 
elements from $1+A_{>0}$ on the diagonal. Hence, up to quotients of 
elements from $1+A_{>0}$, the principal minors of $\Omega$ are products 
of terms $\langle \tilde{f}_t, \tilde{f}_t\rangle_{\lambda} \in 1+A_{>0}$,  
for various~$t$. We conclude that each principal minor of $\Omega$ is a 
quotient of elements from $1+A_{>0}$. Since, on the other hand, the 
coefficients of $\Omega$ lie in $A$, so do all principal minors of $\Omega$. 
Hence (F2) holds. So Lemma~\ref{critnorm} shows that $\bH_n$ is integral. 
\end{exmp}

\begin{rem} \label{typA} In the above setting, consider the special 
case where $\Gamma=\Z$ and $A$ is the ring of Laurent polynomials in
$v=e^1$. The parabolic subgroup $\fS_n:=\langle s_1,\ldots,s_{n-1}
\rangle \subseteq W_n$ can be identified with the symmetric group on 
$\{1,\ldots,n\}$ where $s_i$ corresponds to the basic transposition 
$(i,i+1)$ for $1\leq i \leq n-1$. Let $\bH(\fS_n) \subseteq \bH_n$ be 
the corresponding parabolic subalgebra.  We have the following diagram:
\begin{center}
\begin{picture}(250,50)
\put(  3, 25){$\fS_n$}
\put(  4, 05){$L$ :}
\put( 61, 25){\circle{10}}
\put( 66, 25){\line(1,0){29}}
\put(100, 25){\circle{10}}
\put(105, 25){\line(1,0){20}}
\put(135, 22){$\cdot$}
\put(145, 22){$\cdot$}
\put(155, 22){$\cdot$}
\put(165, 25){\line(1,0){20}}
\put(190, 25){\circle{10}}
\put( 56, 37){$s_1$}
\put( 58, 05){$a$}
\put( 96, 37){$s_2$}
\put( 98, 05){$a$}
\put(183, 37){$s_{n-1}$}
\put(188, 05){$a$}
\end{picture}
\end{center}
Now, by Hoefsmit \cite[\S 2.3]{Hoefs}, we obtain a complete set of 
irreducible representations of $\bH(\fS_n)_K$ by restricting the irreducible 
representations of $\bH_{n,K}$ with label of the form $\lambda=(\lambda_1,
\varnothing)$, where $\lambda_1$ is a partition of $n$. Hence, the
fact that $\bH_n$ is integral and normalized immediately implies that
$\bH(\fS_n)$ is integral and normalized, too. Alternatively, one could
also work directly with the Dipper--James construction of Hoefsmit's 
matrices for type $A_{n-1}$ in \cite[Theorem~4.9]{DiJa2}.
\end{rem}

The following arguments are inspired from the proof of
\cite[Theorem~4.10]{my02}.

\begin{lem} \label{lem01} Assume that $\bH$ is integral and normalized. 
\begin{itemize}
\item[(a)] We have $c_{w,\lambda}^{ij}\in \{0,\pm 1\}$ for all $w\in W$,
$\lambda \in \Lambda$ and $1\leq i,j \leq d_\lambda$.
\item[(b)] For any $\lambda\in\Lambda$ and $1\leq i,j \leq d_\lambda$, 
there exists a unique $w\in W$ such that $c_{w,\lambda}^{ij} \neq 0$;
we denote that element by $w=w_\lambda(i,j)$. The correspondence 
$(\lambda,i,j)\mapsto w_\lambda(i,j)$ defines a bijective map
\[ \{(\lambda,i,j)\} \mid \lambda \in \Lambda,1\leq i,j \leq d_\lambda\}
\longrightarrow W.\]
\end{itemize}
In particular, the sets $\cC_\lambda$ defined above form a partition 
of $W$:
\[ W=\coprod_{\lambda \in \Lambda} \cC_\lambda \qquad \mbox{and}
\qquad |\cC_\lambda|=d_\lambda^2 \quad \mbox{for all $\lambda \in
\Lambda$}.\]
\end{lem}

\begin{proof} First we construct the desired map in (b). Fix $\lambda\in 
\Lambda$ and $1\leq i,j \leq d_\lambda$. Then consider the Schur relation 
where $\lambda=\mu$, $i=l$ and $j=k$:
\[ \sum_{w \in W} (c_{w,\lambda}^{ij})^2=r_\lambda=1.\]
Since the leading matrix coefficients are integers, we conclude that
there exists a unique $w=w_\lambda(i,j)$ such that $c_{w,\lambda}^{ij}=
\pm 1$ and $c_{y,\lambda}^{ij}=0$ for all $y\in W\setminus \{w\}$. Thus,
we have a map $(\lambda,i,j)\mapsto w=w_\lambda(i,j)$; furthermore, note 
that once we have shown that this map is bijective, then (a) follows.

Next we show that the above map is surjective. Let $w\in W$. Then the 
``second Schur relations'' show that there exist some $\lambda \in 
\Lambda$ and $1\leq i,j \leq d_\lambda$ such that $c_{w,\lambda}^{ij}
\neq 0$. The previous argument implies that $w=w_\lambda(i,j)$. Thus, 
the above map is surjective. Since $\dim \bH_K=|W|=\sum_{\lambda\in 
\Lambda} d_\lambda^2$, that map is between finits sets of the same 
cardinality. Hence, the map is bijective.
\end{proof}

\begin{rem} \label{incells} In the setting of Lemma~\ref{lem01}, let
$\lambda \in \Lambda$ and consider the set $\cC_\lambda$. First of all, we 
have a unique labelling of the elements in $\cC_\lambda$:
\[ \cC_{\lambda}=\{w_\lambda(i,j) \mid 1\leq i,j \leq d_\lambda\}.\]
It follows from \cite[Theorem~4.4(b)]{my02} that 
\begin{equation*}
w_\lambda(i,j)^{-1}=w_\lambda(j,i).\tag{a}
\end{equation*}
In particular, $w_\lambda(i,j)$ is an involution if and only of $i=j$.
Furthermore, let $i,j \in \{1,\ldots,d_\lambda\}$ and define
\begin{equation*}
{^j}\cC_\lambda:=\{w_\lambda(k,j)\mid 1\leq k \leq d_\lambda\} \quad
\mbox{and}\quad \cC_\lambda^i:=\{w_\lambda(i,l)\mid 1\leq l \leq d_\lambda\}.
\tag{b}
\end{equation*}
It is shown in \cite[Theorem~4.10]{my02} that ${^j}\cC_\lambda$ is contained 
in a left cell of $W$ and $\cC_\lambda^j$ is contained in a right cell of $W$.
In particular, the whole set $\cC_\lambda$ is contained in a two-sided cell
of $W$. 
\end{rem} 

\begin{lem} \label{lem01a} Assume that $\bH$ is integral and normalized.
Furthermore, assume that $\chi_{\fC} \in \Irr(\bH_K)$ for all left cells
$\fC$ of $W$. Then the sets 
\[ \{{^j}\cC_\lambda \mid \lambda \in \Lambda, 1\leq j \leq d_\lambda\}\]
are precisely the left cells of $W$. The character of the left cell
representation afforded by ${^j}\cC_\lambda$ is given by $\chi_\lambda$.
Furthermore, each left cell contains a unique involution.
\end{lem}

\begin{proof} Let $\lambda\in \Lambda$ and $1\leq j \leq d_\lambda$. By
Remark~\ref{incells}, there is a left cell $\fC$ of $W$ such that 
${^j}\cC_\lambda \subseteq \fC$. Now Proposition~\ref{mylem} shows that 
$\chi_\lambda$ occurs with non-zero multiplicity in $\chi_{\fC}$. Since
$\chi_{\fC} \in \Irr(\bH_K)$, we conclude that $\chi_{\fC}=\chi_\lambda$.
In particular, this means that the underlying $\bH_K$-modules have the
same dimension and so $|\fC|=\chi_\lambda(1)=d_\lambda$. Consequently,
we have $\fC={^j}\cC_\lambda$. Thus, we have shown that each set 
${^j}\cC_\lambda$ is a left cell of $W$. Since these sets form a partition
of $W$, we conclude that they are precisely the left cells of $W$.
The statement concerning involutions now follows from Remark~\ref{incells}(a).
\end{proof}

\begin{exmp} \label{bn1} Let $W=W_n$ and $L\colon W_n\rightarrow \Gamma$ 
be as in Example~\ref{asym} (the ``asymptotic case'' in type $B_n$). We 
have already noted in Example~\ref{bn} that $\bH_n$ is integral and 
normalized. Thus, by Lemma~\ref{lem01}, we have a partition
\[W_n=\coprod_{(\lambda_1,\lambda_2)\in\Lambda_n}\cC_{(\lambda_1,
\lambda_2)}.\]
Let us identify the sets $\cC_{(\lambda_1,\lambda_2)}$. For this purpose, 
we need the results of Bonnaf\'e--Iancu \cite{BI} concerning the left 
cells of $W_n$ and their characters. These results remain valid in the 
``asymptotic case'' by Bonnaf\'e \cite[\S 5]{BI2}. More precisely, in
\cite[Theorem~7.7]{BI}, the left cells in the ``generic asymptotic case''
(that is, with respect to the weight function $L_0\colon W_n \rightarrow 
\Gamma_0$) are described in terms of a generalized Robinson--Schensted 
correspondence.  In \cite[Cor.~5.2]{BI2}, it is shown that two elements of 
$W_n$ lie in the same left cell with respect to $L\colon W_n \rightarrow
\Gamma$ if and only if this is the case with respect to $L_0$. Hence the 
combinatorial description of the left cells remains valid in the 
``asymptotic case''. The fact that the characters afforded by the left 
cells are all irreducible now follows by exactly the same argument as in 
\cite[Prop.~7.9]{BI}. 

Now let $\fC$ be a left cell. Then $\fC$ is precisely the set of elements in 
a generalized Robinson--Schensted cell (RS-cell for short). Such a cell 
is labelled by a pair of bitableaux of the same shape and size~$n$. Let 
$(\lambda_1,\lambda_2)\in \Lambda_n$ be the bipartition specifying the 
shape of the bitableaux. Then we have 
\[ \chi_{\fC}=\chi_{(\lambda_1,\lambda_2^*)}\in \Irr(\bH_{n,K}).\]
(Note that, in \cite[Prop.~7.11]{BI}, the labelling is given by $(\lambda_2,
\lambda_1^*)$. The reason for the different labelling is that, in 
\cite{BI}, the left cell representations are defined using the $C'$-basis 
while here we use the $C$-basis. As explained in Remark~\ref{otherrep}, the 
character of the left cell representation defined with respect to the 
$C'$-basis is given by $\chi_{\fC} \circ \delta$. On the level of characters
of $W_n$, this corresponds to tensoring with the sign character; see 
\cite[9.4.1]{ourbuch} for a precise statement and more details. The 
effect of tensoring with the sign character is described in 
\cite[5.5.6]{ourbuch}.) Now we claim that 
\[ \cC_{(\lambda_1,\lambda_2^*)}=\{w\in W\mid \mbox{ $w$ belongs to an 
RS-cell of shape $(\lambda_1,\lambda_2)$}\}.\]
Indeed, let $w\in \cC_{(\lambda_1,\lambda_2^*)}$. Then, by 
Proposition~\ref{mylem}, the character $\chi_{(\lambda_1,\lambda_2^*)}$ 
occurs with non-zero multiplicity in $\chi_{\fC}$, where $\fC$ is the left 
cell containing $w$. So the above formula for $\chi_{\fC}$ shows that $w$ 
belongs to an RS-cell of shape $(\lambda_1,\lambda_2)$. The reverse 
implication now follows formally from the fact that both the sets 
$\cC_{(\lambda_1,\lambda_2)}$ and the RS-cells form a partition of $W$. 
\end{exmp}

\section{On Lusztig's conjectures} \label{int-nor}

We keep the setting of the previous section. Now our aim is to develop some
tools and criteria for attacking the properties in Conjecture~\ref{Pconj}.
Throughout this section, we assume that
\begin{center}
{\em $\bH$ is integral and normalized (see Definition~\ref{def01})}.
\end{center}
We begin with an alternative characterization of the $\ba$-function. 
First note that
\[ \ba(z)=\min \{\gamma \geq 0 \mid e^\gamma h_{x,y,z}\in \cO \mbox{ for
all $x,y \in W$}\}\]
for any $z \in W$. This simply follows from the equality $\cO\cap A=
A_{\geq 0}$. 

\begin{prop} \label{prop1} For any $z\in W$, we have 
\[{\ba}(z) =\min \{\gamma \geq 0\mid e^\gamma \, \fX_\lambda^{ij}
(D_{z^{-1}}) \in \cO \mbox{ for all $\lambda\in \Lambda$, $1 \leq i,j 
\leq d_\lambda$}\}.\]
Furthermore, if $z\in \cC_{\lambda_0}$ (see Lemma~\ref{lem01}), we have 
$\alpha_{\lambda_0} \leq \ba(z)$.
\end{prop}

\begin{proof} Let $\gamma_0\geq 0$ be minimal such that $e^{\gamma_0}\, 
\fX_\lambda^{ij}(D_{z^{-1}})\in \cO$ for all $\lambda \in \Lambda$, 
$1 \leq i,j \leq d_\lambda$. First we show that $\ba(z)\leq \gamma_0$. 
For this purpose, it is enough to show that $e^{\gamma_0} h_{x,y,z}\in\cO$ 
for all $x,y\in W$. Let $x,y\in W$. To evaluate $h_{x,y,z}$, we use 
the formula
\begin{align*}
\varepsilon_x\,\varepsilon_y\,\varepsilon_z\,h_{x,y,z}&=
\tau(C_xC_yD_{z^{-1}})=\sum_{\lambda \in \Lambda}
\frac{1}{c_\lambda} \, \chi_\lambda(C_xC_yD_{z^{-1}})\\ &=
\sum_{\lambda \in \Lambda} \frac{1}{c_\lambda} \mbox{trace} \Bigl
(\fX_\lambda(C_x)\fX_{\lambda}(C_y)\fX_{\lambda}(D_{z^{-1}})\Bigr)\\&=
\sum_{\lambda \in \Lambda} \sum_{i,j,k=1}^{d_\lambda} \frac{1}{c_\lambda} \,
\fX_\lambda^{ij}(C_x)\,\fX_{\lambda}^{jk}(C_y)\,\fX_{\lambda}^{ki}
(D_{z^{-1}}).
\end{align*}
Now, for any $\lambda \in \Lambda$, we have $r_\lambda=1$ and so 
$c_\lambda=e^{-2\alpha_\lambda} f_\lambda$ where $f_\lambda \in 
{\fI}_{>0}$. Thus, we obtain
\[ e^{\gamma_0}h_{x,y,z}=\pm\sum_{\lambda\in\Lambda}
\sum_{i,j,k=1}^{d_\lambda} f_\lambda^{-1}\Bigl(e^{\alpha_\lambda}
\fX_\lambda^{ij}(C_x)\Bigr) \Bigl(e^{\alpha_\lambda}\, 
\fX_{\lambda}^{jk}(C_y)\Bigr)\Bigl(e^{\gamma_0}\fX_{\lambda}^{ki}
(D_{z^{-1}})\Bigr).\]
As we have already noted, for all $w\in W$, we have 
\[e^{\alpha_\lambda}\,\fX_\lambda^{ij}(C_w) \equiv \varepsilon_w\,
c_{w,\lambda}^{ij} \quad \bmod \fp.\]
Hence, since $f_\lambda^{-1}\in \cO$, all terms on the right hand side of 
the above identity lie in $\cO$. So we obtain $e^{\gamma_0}h_{x,y,z}\in \cO$ 
as desired. Next we show that $\gamma_0\leq \ba(z)$. By the definition of 
$\gamma_0$, there exists some $\lambda_0\in \Lambda$ such that 
\[ e^{\gamma_0} \fX_{\lambda_0}^{k_0i_0} (D_{z^{-1}}) \not\equiv 0
\quad \bmod \fp \qquad \mbox{for some $1\leq i_0,k_0\leq d_{\lambda_0}$}.\]
We also know by Lemma~\ref{lem01} that there exist $x_0,y_0\in W$ such that 
$c_{x_0,\lambda_0}^{i_0k_0}=\pm 1$ and $c_{y_0,\lambda_0}^{k_0k_0}=\pm 1$. 
As before, we obtain an identity
\[ e^{\gamma_0}h_{x_0,y_0,z}=\pm 
\sum_{\lambda\in\Lambda}\sum_{i,j,k=1}^{d_\lambda} 
f_\lambda^{-1}\, \Bigl(e^{\alpha_\lambda}\fX_\lambda^{ij}(C_{x_0})\Bigr)
\Bigl(e^{\alpha_\lambda}\fX_{\lambda}^{jk}(C_{y_0})\Bigr)\Bigl(e^{\gamma_0}
\fX_{\lambda}^{ki}(D_{z^{-1}})\Bigr).\]
All terms on the right hand side lie in $\cO$. So we obtain
\[ e^{\gamma_0}h_{x_0,y_0,z} \equiv \pm \sum_{\lambda \in\Lambda} 
\sum_{i,j,k=1}^{d_\lambda} c_{x_0,\lambda}^{ij}\,c_{y_0,\lambda}^{jk}\, 
\Bigl(e^{\gamma_0}\fX_{\lambda}^{ki} (D_{z^{-1}})\Bigr) \bmod \fp;\]
note that $f_\lambda^{-1}\equiv 1\bmod \fp$ since $f_\lambda \in \fI_{>0}$.
Now, by Lemma~\ref{lem01}, we have 
\begin{align*}
 c_{x_0,\lambda}^{ij} &\neq 0 \Rightarrow (\lambda,i,j)=
(\lambda_0,i_0,k_0),\\
c_{y_0,\lambda}^{jk} & \neq 0 \Rightarrow (\lambda,j,k)=
(\lambda_0,k_0,k_0).
\end{align*}
Hence we obtain 
\begin{align*}
 e^{\gamma_0}h_{x_0,y_0,z} &\equiv \pm c_{x_0,\lambda_0}^{i_0k_0}\,
c_{y_0,\lambda_0}^{k_0k_0}\, \Bigl(e^{\gamma_0}\fX_{\lambda_0}^{k_0i_0} 
(D_{z^{-1}})\Bigr) \\& \equiv \pm e^{\gamma_0}\fX_{\lambda_0}^{k_0i_0} 
(D_{z^{-1}})\not\equiv 0 \bmod \fp
\end{align*}
Consequently, we must have $\gamma_0 \leq \ba(z)$
and so $\gamma_0=\ba(z)$.

Finally, if $z\in \cC_{\lambda_0}$, then $c_{z,\lambda_0}^{ij}=\pm 1$ for 
some $i,j$. Since 
\[ \varepsilon_z\,e^{\alpha_{\lambda_0}} \fX_{\lambda_0}^{ji}(D_{z^{-1}}) 
\equiv c_{z^{-1}, \lambda_0}^{ji} \equiv c_{z,\lambda_0}^{ij} \equiv 
\pm 1\quad \bmod \fp,\]
we conclude that $\ba(z)\geq \alpha_{\lambda_0}$. 
\end{proof}

\begin{defn} \label{def2} Recall that $\bH$ is assumed to be normalized
and integral. Let $z\in W$. Then we set 
\[ \alpha_z:=\alpha_\lambda,\] 
where $\lambda$ is the unique element of $\Lambda$ such that $z \in
\cC_\lambda$; see Lemma~\ref{lem01}. By Proposition~\ref{prop1}, we have
$\alpha_z\leq \ba(z)$ for all $z\in W$.
\end{defn}

The following result shows that the identity $\ba(z)=\alpha_z$ holds for
$z\in W$ once we know that (P4) from the list of Lusztig's conjectures holds.
Note, however, that it seems to be very hard to prove (P4) directly. 
Therefore, in Lemma~\ref{cor1a} below, we shall establish a somewhat 
different criterion for proving that the identity $\ba(z)=\alpha_z$ holds 
for all $z\in W$ (and it is this latter criterion which will be used in 
Section~5 in dealing with type $B_n$). 

\begin{lem} \label{cor1b} Assume that (P4) holds. Then $\ba(z)=\alpha_z$ 
for all $z\in W$.
\end{lem}

\begin{proof} Let $z\in W$. We already know by Proposition~\ref{prop1} 
that $\alpha_z\leq \ba(z)$. To prove the reverse inequality, let $x,y\in W$ 
be such that $e^{\ba(z)}\, h_{x,y,z}\in A_{\geq 0}$ has a non-zero 
constant term. As in the proof of Proposition~\ref{prop1}, we have 
\[ e^{\ba(z)}h_{x,y,z} \equiv \pm \sum_{\lambda \in 
\Lambda} \sum_{i,j,k=1}^{d_\lambda} c_{x,\lambda}^{ij} \,
c_{y,\lambda}^{jk}\, \Bigl(e^{\ba(z)}\,\fX_{\lambda}^{ki} 
(D_{z^{-1}})\Bigr) \bmod \fp.\]
We are assuming that the left hand side is $\not\equiv 0 \bmod \fp$. So 
there exists some $\lambda_0\in \Lambda$ and $1\leq i,j,k\leq d_{\lambda_0}$ 
such that $c_{x,\lambda_0}^{ij} \neq 0$, $c_{y,\lambda_0}^{jk}\neq 0$ and 
\begin{equation*}
e^{\ba(z)}\,\fX_{\lambda_0}^{ki}(D_{z^{-1}}) \not\equiv 0 \bmod\fp.
\tag{$\dagger$}
\end{equation*}
The condition ($\dagger$) implies that $\alpha_{\lambda_0} \geq \ba(z)$. 
Furthermore, we have $x,y \in \cC_{\lambda_0}$ and so $\ba(x)\geq 
\alpha_{\lambda_0}$, $\ba(y)\geq \alpha_{\lambda_0}$ by 
Proposition~\ref{prop1}.  Hence we obtain 
\[ \ba(x)\geq \alpha_{\lambda_0} \geq \ba (z).\]
But, since $h_{x,y,z}\neq 0$, we have $z \leq_{\cR} x$ and so $\ba(x)\leq
\ba(z)$, thanks to the assumption that (P4) holds. Thus, we conclude
that $\ba(z)=\alpha_{\lambda_0}$. But then ($\dagger$) also yields that 
$c_{z,\lambda_0}^{ik}=c_{z^{-1},\lambda_0}^{ki} \neq 0$ and so 
$z\in \cC_{\lambda_0}$. Thus, we have $\ba(z)=\alpha_{\lambda_0}=\alpha_z$, 
as claimed.
\end{proof}

\begin{lem} \label{cor1a} Assume that the following implication holds
for any $x,y \in W$:
\begin{equation*}
x \leq_{\cLR} y \qquad \Rightarrow \qquad \alpha_y \leq \alpha_x.
\tag{$\clubsuit$}
\end{equation*}
Then we have $\ba(z)=\alpha_z$ for all $z\in W$ (and, consequently, 
(P4) holds).
\end{lem}

\begin{proof} Let $z\in W$. By Proposition~\ref{prop1}, we already know 
that $\alpha_z\leq \ba(z)$; furthermore, in order to prove equality, it will
be enough to show that 
\[ e^{\alpha_z}\fX_\lambda^{ij}(D_{z^{-1}})\in \cO \qquad \mbox{for 
all $\lambda \in \Lambda$ and $1\leq i,j \leq d_\lambda$}.\]
To prove this, let $\lambda\in \Lambda$ be such that $\fX_\lambda^{ij}
(D_{z^{-1}}) \neq 0$ for some $i,j$.  Let $\fC$ be a left cell such that
$[\chi_{\fC}:\chi_\lambda] \neq 0$, that is, $\chi_\lambda$ occurs with 
non-zero multiplicity in the character afforded by $\fC$. Then 
$\fX_\lambda$ will occur (up to equivalence) in the decomposition of 
$\fX_{\fC}$ as a direct sum of irreducible representations. Hence, since 
$\fX_\lambda(D_{z^{-1}})\neq 0$, we will also have $\fX_{\fC}(D_{z^{-1}}) 
\neq 0$. Recalling the definition of $\fX_{\fC}$, we deduce that there 
exists some $x \in \fC$ such that $D_{z^{-1}}C_x \neq 0$. Since $\tau$ is 
non-degenerate, we have 
\[ \tau(D_{z^{-1}}C_xC_w)\neq 0 \qquad \mbox{for some $w\in W$}.\]
This yields $\pm h_{x,w,z}=\tau(C_xC_wD_{z^{-1}})=\tau(D_{z^{-1}}C_xC_w)
\neq 0$ and so $z \leq_{\cR} x$. Hence ($\clubsuit)$ implies that $\alpha_x 
\leq \alpha_z$. Now we claim that $\alpha_x= \alpha_\lambda$. This can be 
seen as follows. Let $i \in \{1,\ldots,d_\lambda\}$. Then the right hand 
side of the formula in Proposition~\ref{mylem} is non-zero and so there 
exists some $y \in \fC$ and some $k$ such that $c_{y, \lambda}^{ik} \neq 0$. 
Consequently, we have $y \in \cC_\lambda$ and so $\alpha_y= \alpha_\lambda$. 
On the other hand, $(\clubsuit$) implies that the function $w \mapsto 
\alpha_w$ is constant on two-sided cells. Hence we can deduce that 
$\alpha_x=\alpha_y= \alpha_\lambda$, as claimed.

Now the inequality $\alpha_x \leq  \alpha_z$ yields $\alpha_\lambda \leq 
\alpha_z$ and so 
\[ e^{\alpha_z}\fX_\lambda^{ij}(D_{z^{-1}})= e^{\alpha_\lambda}
\fX_\lambda^{ij}(D_{z^{-1}}) \, e^{\alpha_z-\alpha_\lambda}\in \cO,\]
as desired.
\end{proof}

We will now show that the validity of the identity $\ba(z)=\alpha_z$ for 
all $z \in W$ together with {(P4)} formally imply most of the
remaining conjectures from the list {(P1)--(P15)}.

We introduce the following notation. Let $x,y\in W$. By Lemma~\ref{lem01},
we can write $x=w_\lambda(i_0,j_0)$ and $y=w_\mu(k_0,l_0)$ where $\lambda,
\mu\in \Lambda$, $1\leq i_0,j_0 \leq d_\lambda$ and $1\leq k_0,l_0 \leq 
d_\mu$ are uniquely determined. Then we set 
\[ z=x\star y:=w_\lambda(l_0,i_0)  \qquad \mbox{if $\lambda=\mu$ and $j_0=
k_0$}.\]
We note the following equivalences for $x,y,z\in \cC_\lambda$:
\begin{equation*}
 z=x\star  y \quad \Leftrightarrow \quad x=y\star z \quad 
\Leftrightarrow \quad y=z \star x.\tag{$\heartsuit$}
\end{equation*}
If $z=x\star y$ as above, we also set 
\[ n_{x,y,z}=c_{x,\lambda}^{i_0j_0} c_{y,\lambda}^{j_0l_0}
c_{z,\lambda}^{l_0i_0}=\pm 1.\] 
Note that ($\heartsuit$) implies $n_{x,y,z}=n_{y,z,x}=n_{z,x,y}$.
With this notation, we now have the following result.

\begin{lem} \label{cor1b1} Assume that (P4) holds. Let $x,y,z\in W$. 
If there is no $\lambda \in \Lambda$ such that $x,y,z \in \cC_\lambda$,
then $\gamma_{x,y,z}=0$. If $x,y,z\in \cC_\lambda$, then 
\[ \gamma_{x,y,z}=\gamma_{y,z,x}=\gamma_{z,x,y}=\left\{\begin{array}{cl} 
n_{x,y,z}&\qquad \mbox{if $z=x\star y$},\\ 0 & \qquad 
\mbox{otherwise}.\end{array}\right.\]
\end{lem}

\begin{proof} By Remark~\ref{invers}, we have $\ba(z)=\ba(z^{-1})$. 
This yields
\[\gamma_{x,y,z} \equiv e^{\ba(z)}h_{x,y,z^{-1}} \quad \bmod \fp.\]
Writing $h_{x,y,z^{-1}}=\varepsilon_x\,\varepsilon_y\,\varepsilon_z\,
\tau(C_xC_yD_z)$ and $\tau=\sum_\mu c_\mu^{-1} \chi_\mu$, we obtain 
\[\gamma_{x,y,z} \equiv \sum_{\mu\in \Lambda} \sum_{i,j,k=1}^{d_\mu} 
c_{x,\mu}^{ij} \, c_{y,\mu}^{jk}\, \Bigl(\varepsilon_z\,e^{\ba(z)}\, 
\fX_{\mu}^{ki}(D_z)\Bigr) \quad \bmod \fp;\]
see the argument in the proof of  Proposition~\ref{prop1}. Now, if there is 
no $\lambda$ such that $x,y\in \cC_\lambda$, then the above sum certainly 
is zero. 

Hence it remains to consider the case where $x,y\in \cC_\lambda$ for some 
$\lambda \in \Lambda$ (which is uniquely determined). We write
$x=w_\lambda(i_0,j_0)$ and $y=w_\lambda(k_0,l_0)$, where $1\leq i_0,j_0,
k_0,l_0\leq d_\lambda$ are uniquely determined. Then the above sum reduces to 
\[\gamma_{x,y,z} \equiv \delta_{j_0k_0} \,c_{x,\lambda}^{i_0j_0}\, 
c_{y,\lambda}^{k_0l_0}\, \Bigl(\varepsilon_z\,e^{\ba(z)}\, 
\fX_{\lambda}^{l_0i_0}(D_z) \Bigr) \quad \bmod \fp.\]
First assume that $z=x\star y$, that is, we have $j_0=k_0$ and $z=
w_\lambda(l_0,i_0)$. In particular, $z\in \cC_\lambda$ and so $\ba(z)=
\alpha_\lambda$ by Lemma~\ref{cor1b}. Hence, we obtain
\[ \gamma_{x,y,z} =c_{x,\lambda}^{i_0j_0}\, c_{y,\lambda}^{k_0l_0}\, 
c_{z,\lambda}^{l_0i_0}=n_{x,y,z}=\pm 1,\]
as desired. 

Conversely, assume that $\gamma_{x,y,z}\neq 0$. First of all, this means 
that $j_0=k_0$. Furthermore, the constant term of $e^{\ba(z)}\, 
\fX_{\lambda}^{l_0i_0}(D_z)$ is non-zero and so $\alpha_{\lambda} \geq 
\ba(z)$. On the other hand, we have $h_{x,y,z^{-1}} \neq 0$ and so 
$z^{-1} \leq_{\cR} x$. Hence (P4) yields $\ba(z)=\ba(z^{-1})\geq \ba(x)=
\alpha_\lambda$ (since $x\in \cC_\lambda$) and so $\ba(z)=\alpha_\lambda$. 
Consequently, the constant term of $\varepsilon_z\,e^{\ba(z)}\, 
\fX_{\lambda}^{l_0i_0} (D_z)$ is $c_{z,\lambda}^{l_0i_0}$ and so
\[0 \neq \gamma_{x,y,z}=\,c_{x,\lambda}^{i_0j_0}\,
c_{y,\lambda}^{j_0l_0}\, c_{z,\lambda}^{l_0i_0}.\] 
Hence, we must have $j_0=k_0$ and $z=x\star y$.

The identity $\gamma_{x,y,z}=\gamma_{y,z,x}=\gamma_{z,x,y}$ 
easily follows from the symmetry in ($\heartsuit$).
\end{proof}

\begin{lem} \label{cor1c} Assume that (P4) holds.
Then (P1), (P2), (P3), (P5), (P6), (P7), (P8) and  (P14) hold. 
\end{lem}

\begin{proof}  {\bf (P1)} Let $z\in W$. We must show $\ba(z)\leq \bD(z)$. 
Now, by the definition of the symmetrizing trace  $\tau$, we have 
\[  \tau(C_z)=\varepsilon_z\,\overline{P}_{1,z}^*.\]
On the other hand, we have $\tau=\sum_\lambda c_\lambda^{-1}\chi_\lambda$.
Using the expression $c_\lambda=e^{-2\alpha_\lambda}f_\lambda$ where
$f_\lambda^{-1}\in \cO$, this yields the identity
\[ \varepsilon_z\,\overline{P}_{1,z}^*=\sum_{\lambda\in \Lambda} 
\sum_{i=1}^{d_\lambda} f_\lambda^{-1} (e^{\alpha_\lambda} \, 
\fX_\lambda^{ii}(C_z))\, e^{\alpha_\lambda}.\]
Assume that the term in the sum corresponding to $\lambda \in \Lambda$ and 
$1\leq i \leq d_\lambda$ is non-zero. Let $\fC$ be a left cell such that
$[\chi_{\fC}:\chi_\lambda] \neq 0$, that is, $\chi_\lambda$ occurs with
non-zero multiplicity in the character afforded by $\fC$. Then $\fX_\lambda$ 
will occur (up to equivalence) in the decomposition of $\fX_{\fC}$ as a 
direct sum of irreducible representations. Hence, since $\fX_\lambda(C_z)
\neq 0$, we will also have $\fX_{\fC}(C_z) \neq 0$. Recalling the definition 
of $\fX_{\fC}$, we deduce that there exist some $x,y \in \fC$ such that 
$h_{z,x,y}\neq 0$. But then we have $y \leq_{\cR} z$ and so $\ba(z)\leq 
\ba(y)$, by (P4). Now, as in the proof of Lemma~\ref{cor1a}, we 
conclude that $\ba(y)=\ba(y')=\alpha_{\lambda}$, where $y' \in \fC$ is
chosen such that $c_{y',\lambda}^{ik}\neq 0$ for some $k$. Hence, we have 
$\ba(z)\leq \alpha_\lambda$ for all non-zero terms in the above sum. 
Thus, we obtain
\[ \varepsilon_z\,e^{-\ba(z)}\,\overline{P}_{1,z}^*=\sum_{\atop{\lambda\in 
\Lambda}{\ba(z)\leq \alpha_\lambda}} \sum_{i=1}^{d_\lambda} 
f_\lambda^{-1} (e^{\alpha_\lambda} \, \fX_\lambda^{ii}(C_z))\, 
e^{\alpha_\lambda-\ba(z)} \in \cO.\]
Since the left hand side lies in $A$, we conclude that $e^{-\ba(z)}\,
\overline{P}_{1,z}^* \in A_{\geq 0}$ and so $\ba(z) \leq \bD(z)$, as desired.

\medskip
{\bf (P2)}, {\bf (P5)} Let $x,y\in W$ and $d\in \cD$ be such that 
$\gamma_{x,y,d}\neq 0$. We must show that $x=y^{-1}$ and $\gamma_{y^{-1},
y,d} =n_d=\pm 1$. Now, using the expression of $C_w$ in terms of the 
$T$-basis, we obtain 
\begin{align*}
\tau(C_xC_y)&=\sum_{z\in W} \varepsilon_x\,\varepsilon_y\,\varepsilon_z\,
h_{x,y,z}\tau(C_z)=\sum_{z \in W} \varepsilon_x\,\varepsilon_y\,
h_{x,y,z}\, \overline{P}_{1,z}^*\\ &=\sum_{z \in W} \varepsilon_x\,
\varepsilon_y\,e^{\ba(z)}\,h_{x,y,z}\, (e^{-\ba(z)}\,\overline{P}_{1,z}^*).
\end{align*}
By the definition of $\ba(z)$ and (P1) (see the proof above), we have 
$e^{\ba(z)}\,h_{x,y,z}\in A_{\geq 0}$ and $e^{-\ba(z)}\,
\overline{P}_{1,z}^* \in A_{\geq 0}$. So the terms of the above sum 
lie in $A_{\geq 0}$ and we have 
\[ \tau(C_xC_y) \equiv \sum_{z \in \cD} \varepsilon_x\,\varepsilon_y\,
\gamma_{x,y,z^{-1}}\, n_z \quad \bmod \fp.\]
By Remark~\ref{invers}, we have $n_z=n_{z^{-1}}$ and $\cD^{-1}=\cD$.
So the above congruence can also be written in the form 
\[ \tau(C_xC_y) \equiv \sum_{z \in \cD} \varepsilon_x\,\varepsilon_y\,
\gamma_{x,y,z}\, n_z \quad \bmod \fp.\]
Now, Lemma~\ref{cor1b1} shows that the only non-zero term in the 
above sum is $\gamma_{x,y,d}$ and that $d=x\star y$. So we have 
\[ \tau(C_xC_y) \equiv \varepsilon_x\,\varepsilon_y\,
\gamma_{x,y,d}\, n_d \not\equiv 0 \quad \bmod \fp.\]
On the other hand, we also have 
\[ \tau(C_xC_y) \equiv \delta_{x^{-1}y} \bmod A_{>0},\]
where $\delta_{x^{-1}y}$ denotes the Kronecker symbol.  (This easily follows 
from the defining formulas; see \cite[14.5(a)]{Lusztig03}.)  Hence we
obtain the congruence
\[ \delta_{x^{-1}y} \equiv \tau(C_xC_y) \equiv \varepsilon_x\,\varepsilon_y\,
\gamma_{x,y,d}\, n_d \not\equiv 0 \quad \bmod \fp.\]
So we must have $x^{-1}=y$.
But then we also get $\gamma_{y^{-1},y,d}n_d=1$, as required.

\medskip
{\bf (P3)} Let $y \in W$. We want to show that there exists a unique 
$z\in \cD$ such that $\gamma_{y^{-1},y,z}\neq 0$. As in the proof of (P2),
we have 
\[\sum_{z\in \cD} \gamma_{y^{-1},y,z} n_z \equiv \tau(C_{y^{-1}}C_y)
\equiv 1 \quad \bmod \fp.\]
Consequently, there exists some $z\in \cD$ such that $\gamma_{y^{-1},y, z} 
\neq 0$. By Lemma~\ref{cor1b1}, $z=y^{-1} \star y$ is uniquely determined 
with this property.

\medskip 
{\bf (P6)} This is a formal consequence of (P2) and (P3); see 
\cite[14.6]{Lusztig03}. 

\medskip
{\bf (P7)}, {\bf (P8)} This is clear by Lemma~\ref{cor1b1} and 
Remark~\ref{incells}.

\medskip 
{\bf (P14)} This is clear by Remark~\ref{incells}.
\end{proof}

\begin{cor} \label{detnz} We keep the hypotheses of Lemma~\ref{cor1c}.
Let $z\in \cD$ and $\lambda \in \Lambda$ be such that $z\in \cC_\lambda$. 
Then the constant $n_z=\pm 1$ is determined by the formula
\[\varepsilon_z\,e^{\alpha_\lambda}\,\chi_\lambda(T_z) \equiv n_z 
\bmod \fp .\]
Thus, $n_z$ is precisely the leading coefficient of a character value
as defined by Lusztig \cite{Lusztig87}.
\end{cor}

\begin{proof} Since $z \in \cD$, we have $z^2=1$ and so 
$z=w_\lambda(i_0,i_0)$ for a unique $i_0\in \{1,\ldots,d_\lambda\}$. 
Hence we obtain
\[ \varepsilon_z\,e^{\alpha_\lambda} \chi_\lambda(T_z)\equiv \varepsilon_z
\sum_{i=1}^{d_\lambda} e^{\alpha_\lambda}\fX_\lambda^{ii}(T_z)
\equiv \sum_{i=1}^{d_\lambda} c_{z,\lambda}^{ii} \equiv 
c_{z,\lambda}^{i_0i_0} \quad \bmod \fp.\]
On the other hand, by {(P5)} and Lemma~\ref{cor1b1}, we have 
$n_z=\gamma_{z,z,z}=c_{z,\lambda}^{i_0i_0}$, as required.
\end{proof}

\begin{lem} \label{cor1d} We keep the hypotheses of Lemma~\ref{cor1c}
and assume, in addition, that $\chi_{\fC}\in \Irr(\bH_K)$ for every left 
cell $\fC$ of $W$. Then (P13) holds and we have 
\[\cD=\{z\in W\mid z^2=1\}.\] 
\end{lem}

\begin{proof} Let $\fC$ be a left cell of $W$. Let $x\in \fC$. By 
Lemma~\ref{cor1c}, (P3) holds and so there exists a unique $d\in \cD$ 
such that $\gamma_{x^{-1},x,d}\neq 0$. By (P8), we have $d\sim_{\cL} x$ 
and so $d\in \fC$. Thus, each left cell contains an element of $\cD$. 
Note that the above argument also shows that (P13) holds, once we know
that each left cell contains a unique element of $\cD$.

Now, by Lemma~\ref{lem01a}, the total number of left cells equals 
$\sum_{\lambda \in \Lambda} d_\lambda$.  Thus, we have 
\[ |\cD| \geq \sum_{\lambda \in \Lambda} d_\lambda.\]
On the other hand, by a well-known result due to Frobenius--Schur,
the number on the right hand side is the number of all $z\in W$ such
that $z^2=1$. (We also use the fact that every irreducible character of
$W$ can be realized over $\R$; see \cite[6.3.8]{ourbuch}.) Since
$d^2=1$ for all $d\in \cD$ by (P6), we conclude that 
\[ \cD=\{z\in W \mid z^2=1\}\]
and that $\fC$ contains a unique element from $\cD$. 
\end{proof}

Finally, let $J$ be the free abelian group with basis $\{t_w\mid w\in W\}$.
We define a bilinear pairing on $J$ by 
\[ t_x \cdot t_y=\sum_{z\in W} \gamma_{x,y,z^{-1}} \, t_z\qquad
\mbox{for all $x,y \in W$},\]
where the constants $\gamma_{x,y,z^{-1}}\in \Z$ were introduced in
Section~2.

\begin{prop} \label{strucJ} Assume that (P4) holds. Then $J$ is an 
associative ring with identity $1_J=\sum_{z\in \cD} n_zt_z$. For any 
$\lambda \in \Lambda$, we set
\[ J_\lambda:=\langle t_w \mid w \in \cC_\lambda \rangle\subseteq J.\]
Then $J_\lambda$ is a two-sided ideal which is isomorphic to the matrix ring
$M_{d_\lambda}(\Z)$, and we have $J=\bigoplus_{\lambda \in \Lambda}
J_\lambda.$\\ We have $t_z^2=n_zt_z$ for any $z \in \cD$.
\end{prop}

In Example~\ref{expb2}, we will see that negative coefficients actually 
do occur in~$1_J$.

\begin{proof} By Lemma~\ref{cor1c}, we know that {(P1)--(P8)} hold. Hence 
$J$ can be constructed as explained in \cite[Chap.~18]{Lusztig03}. Now 
Lemma~\ref{cor1b1} immediately shows that $J_\lambda$ is a two-sided ideal 
and so we have $J=\bigoplus_\lambda J_\lambda$. To establish the stated 
isomorphism $J_\lambda \cong M_{d_\lambda}(\Z)$, we explicitly construct a 
set of ``matrix units'' in $J_\lambda$. This is done as follows. Let us fix 
$\lambda\in \Lambda$.  For any $1\leq i,j\leq d_\lambda$, we set 
\[ E_\lambda^{ij}:=c_{w,\lambda}^{ij}t_w \qquad
\mbox{where $w=w_\lambda(i,j)$ and $c_{w,\lambda}^{ij} =\pm 1$}.\]
Now let $1\leq i,j,k,l \leq d_\lambda$ and write $x=w_\lambda(i,j)$,
$y=w_\lambda(k,l)$. Then we have 
\[E_\lambda^{ij} \cdot E_\lambda^{kl}= c_{x,\lambda}^{ij}
c_{y,\lambda}^{kl}\,t_x\cdot t_y=c_{x,\lambda}^{ij}
c_{y,\lambda}^{kl}\sum_{z\in W} \gamma_{x,y,z}\,t_{z^{-1}}.  \]
Now Lemma~\ref{cor1b1} shows that the result will be zero unless
$j=k$. So let us finally assume that $j=k$. Then we obtain 
\begin{align*}
E_\lambda^{ij} \cdot E_\lambda^{jl}&= c_{x,\lambda}^{ij} 
c_{y,\lambda}^{jl} \gamma_{x,y,z_0}\, t_{z_0^{-1}},\\ &= 
c_{x,\lambda}^{ij} c_{y,\lambda}^{jl} c_{z_0,\lambda}^{li} 
\gamma_{x,y,z_0} E_\lambda^{il},
\end{align*}
where $z_0:=x\star y=w_\lambda(l,i)$.  By Lemma~\ref{cor1b1}, the 
coefficient of $E_\lambda^{il}$ in the above expression equals~$1$. Thus, 
we have shown that
\[ E_\lambda^{ij}\cdot E_\lambda^{kl}=\delta_{jk}E_\lambda^{il}
\quad \mbox{for $1 \leq i,j,k,l\leq d_\lambda$}.\]
Hence, the elements $E_\lambda^{ij}$ multiply in exactly the same way
as the matrix units in $M_{d_\lambda}(\Z)$, which yields the desired
isomorphism. 

The formula for $t_z^2$, where $z \in \cD$, is obtained as follows. 
We have $t_z^2=\gamma_{z,z,x} t_{x^{-1}}$ where $x=z \star z$ and
$\gamma_{z,z,x}\neq 0$. By (P2), (P6), (P7), we have $\gamma_{x,z,z}
\neq 0$ and so $x=z^{-1}=z$. Then (P5) yields $\gamma_{z,z,z}=n_z= \pm 1$.
\end{proof}

\begin{rem} \label{strucJ1} We keep the setting of 
Proposition~\ref{strucJ}. Let us also assume that $\chi_{\fC} \in 
\Irr(\bH_K)$ for all left cells $\fC$ in $W$. Let $x,y,z\in W$
be such that $x,y,z\in \cC_\lambda$ for some $\lambda$. Then note that, 
by Remark~\ref{incells} and Lemma~\ref{lem01a}, the condition $z=x \star y$ 
is equivalent to the conditions $x \sim_{\cL} y^{-1}$, $y \sim_{\cL} z^{-1}$, 
$z \sim_{\cL} x^{-1}$.  Thus, the multiplication rule in $J$ can now
be formulated as follows:
\[ t_x\cdot t_y=\left\{\begin{array}{cl} \pm t_{z^{-1}} & \qquad \mbox{if 
$x\sim_{\cL} y^{-1}$, $y\sim_{\cL} z^{-1}$, $z\sim_{\cL} x^{-1}$},\\
0 & \qquad \mbox{otherwise}.\end{array}\right.\]
\end{rem}

Summarizing the results in this section, we see that for a normalized
and integral Iwahori--Hecke algebra, property (P4) (or the variant in
Lemma~\ref{cor1a}) implies all the remaining properties in the list of
Lusztig's conjectures except (P9)--(P12) and (P15).

\section{The $\ba$-function in the ``asymptotic case'' in type $B_n$} 
\label{sec:bn}

Throughout this section, we place ourselves in the setting of 
Example~\ref{asym}, where $W_n$ is a Coxeter group of type $B_n$ and
the weight function $L \colon W_n \rightarrow \Gamma$ is given by the
following diagram:
\begin{center}
\begin{picture}(330,50)
\put(  3, 25){$B_n$}
\put(  0, 05){$\{v_s\}$:}
\put( 40, 25){\circle{10}}
\put( 44, 22){\line(1,0){33}}
\put( 44, 28){\line(1,0){33}}
\put( 81, 25){\circle{10}}
\put( 86, 25){\line(1,0){29}}
\put(120, 25){\circle{10}}
\put(125, 25){\line(1,0){20}}
\put(155, 22){$\cdot$}
\put(165, 22){$\cdot$}
\put(175, 22){$\cdot$}
\put(185, 25){\line(1,0){20}}
\put(210, 25){\circle{10}}
\put( 37, 37){$t$}
\put( 36, 05){$e^b$}
\put( 76, 37){$s_1$}
\put( 78, 05){$e^a$}
\put(116, 37){$s_2$}
\put(118, 05){$e^a$}
\put(200, 37){$s_{n-1}$}
\put(208, 05){$e^a$}
\put(240, 20){$b>(n-1)a>0$}
\end{picture}
\end{center}
Let $\bH_n$ be the corresponding Iwahori--Hecke algebra over $A$ where 
$V=e^b$ is the parameter associated with the generator $t$ and $v=e^a$ 
is the parameter associated with the generators $s_1,\ldots,s_{n-1}$ 
of $W_n$. Our aim is to see that we can apply the methods in Section~4. 
In Corollary~\ref{twocor1}, we will be able to show that the key condition
($\clubsuit$) in Lemma~\ref{cor1a} holds in the present setting. 

We shall need some notation from \cite{BI}. Given $w\in W_n$, we denote by 
$\ell_t(w)$ the number of occurrences of the generator $t$ in a reduced 
expression for $w$, and call this the ``$t$-length'' of $w$. 

The parabolic subgroup $\fS_n=\langle s_1,\ldots,s_{n-1}\rangle$ is
naturally isomorphic to the symmetric group on $\{1,\ldots,n\}$, where
$s_i$ corresponds to the basic transposition $(i,i+1)$. For $1 \leq l
\leq n-1$, we set $\Sigma_{l,n-l}:=\{s_1,\ldots,s_{n-1}\}\setminus
\{s_l\}$.  For $l=0$ or $l=n$, we also set $\Sigma_{0,n}=\Sigma_{n,0}=
\{s_1,\ldots,s_{n-1}\}$. Let $Y_{l,n-l}$ be the set of distinguished
left coset representatives of the Young subgroup $\fS_{l,n-l}:=\langle
\Sigma_{l,n-l}\rangle$ in $\fS_n$. We have the parabolic subalgebra
$\bH_{l,n-l}=\langle T_\sigma \mid \sigma \in \fS_{l,n-l}\rangle_A
\subseteq \bH_n$.

We denote by $\leq_{\cL,l}$ the Kazhdan--Lusztig (left) pre-order relation
on $\fS_{l,n-l}$ and by $\sim_{\cL,l}$ the corresponding equivalence
relation. The symbols $\leq_{\cR,l}$, $\leq_{\cLR,l}$, $\sim_{\cR,l}$
and $\sim_{\cLR,l}$ have a similar meaning.

Furthermore, as in \cite[\S 4]{BI}, we set $a_0=1$ and
\[ a_l:=t(s_1t)(s_2s_1t) \cdots (s_{l-1}s_{l-2} \cdots s_1t) \qquad
\mbox{for $l>0$}.\]
Then, by \cite[Prop.~4.4]{BI}, the set $Y_{l,n-l}a_l$ precisely is the 
set of distinguished left coset representatives of $\fS_n$ in $W_n$ whose 
$t$-length equals $l$. Furthermore, every element $w\in W_n$ has a unique 
decomposition
\[w=a_wa_l\sigma_w b_w^{-1} \qquad \mbox{where $l=\ell_t(w)$, $\sigma_w \in
\fS_{l,n-l}$ and $a_w,b_w\in Y_{l,n-l}$};\]
see \cite[4.6]{BI}. With this notation, we have the following result.

\begin{thm}[Bonnaf\'e--Iancu \protect{\cite{BI}} and Bonnaf\'e 
\protect{\cite[\S 5]{BI2}}] \label{mainbi} In the above setting, let 
$x,y\in W_n$. Then we have $x\sim_{\cL} y$ if and only if $l:=\ell_t(x)=
\ell_t(y)$, $b_x=b_y$ and $\sigma_x \sim_{\cL,l} \sigma_y$. 
\end{thm}

(In \cite[Theorem~7.7]{BI}, the above statement is proved in the ``generic 
asymptotic case''. As already discussed in Example~\ref{bn1}, this 
remains valid in the ``asymptotic case'' by \cite[Cor.~5.2]{BI2}.)

We shall also need the following result on the elementary steps in the 
relation $\leq_{\cL}$.

\begin{prop}[Bonnaf\'e--Iancu \protect{\cite{BI}} and Bonnaf\'e 
\protect{\cite[\S 5]{BI2}}] \label{mainbiI} In the above setting, let 
$x,y\in W_n$ be such that $x \leftarrow_{\cL} y$. Then we have $\ell_t(x)=
\ell_t(y)$ or $x=ty>y$. In particular, we always have $\ell_t(y)\leq
\ell_t(x)$.  (A similar result holds for $\leftarrow_{\cR}$.)
\end{prop}

(The precise references are Theorems~6.3, 6.6 and Corollary~6.7 in 
\cite{BI} for the ``generic asymptotic case''. In \cite[Cor.~5.2]{BI2}, 
it is shown that, if $x\leftarrow_{\cL} y$ with respect to $L\colon 
W_n \rightarrow \Gamma$, then we also have $x \leftarrow_{\cL} y$
with respect to $L_0\colon W_n \rightarrow \Gamma_0$ as in 
Example~\ref{asym}. Hence the assertions hold in the 
``asymptotic case'' too.)

As discussed in Example~\ref{bn1}, let 
$\Lambda_n$ be the set of bipartitions of $n$. Then the partition 
\[ W_n=\coprod_{(\lambda_1,\lambda_2)\in \Lambda_n} 
\cC_{(\lambda_1,\lambda_2)}\]
is explicitly given by  
\[\cC_{(\lambda_1,\lambda_2^*)}=\{ w\in W_n \mid \mbox{ $w$ belongs to an
RS-cell of shape $(\lambda_1,\lambda_2)$}\}. \]
Furthermore, for $w\in \cC_{(\lambda_1,\lambda_2^*)}$, we have
\[ \alpha_w=\alpha_{(\lambda_1,\lambda_2^*)}=b\,|\lambda_2|+a\, 
(n(\lambda_1)+ 2n(\lambda_2^*)-n(\lambda_2)).\] 
Finally, we need the following result concerning the relation $\leq_{\cLR}$.

\begin{thm}[Bonnaf\'e \protect{\cite{BI2}}] \label{mainbi2} In the
above setting, let $x,y\in W_n$ be such that $l:=\ell_t(x)=\ell_t(y)$. Then
we have $x\leq_{\cLR} y$ if and only if $\sigma_x \leq_{\cLR,l} \sigma_y$.
Furthermore, the sets $\cC_{(\lambda_1,\lambda_2)}$, $(\lambda_1,\lambda_2)
\in \Lambda_n$, are precisely the two-sided cells of $W_n$.
\end{thm}

Using known results on the two-sided cells in the symmetric group and the 
Robinson--Schensted correspondence, we can translate the above statement
into a combinatorial description of the relation $\leq_{\cLR}$ for $W_n$. 
To state this, we need to introduce some notation. Recall the definition 
of the dominance order on partitions in Example~\ref{symgroup}. Following 
Dipper--James--Murphy \cite[\S 3]{DJM}, we can extend this partial order 
to bipartitions, as follows.

Let $\lambda=(\lambda_1, \lambda_2)$ and $\mu=(\mu_1,\mu_2)$ be bipartitions 
of $n$, with parts 
\begin{gather*}
\lambda_1=(\lambda_1^{(1)} \geq \lambda_1^{(2)} \geq \cdots \geq 0),\qquad
\lambda_2=(\lambda_2^{(1)} \geq \lambda_2^{(2)} \geq \cdots \geq 0),\\
\mu_1=(\mu_1^{(1)} \geq \mu_1^{(2)} \geq \cdots \geq 0),\qquad 
\mu_2=(\mu_2^{(1)} \geq \mu_2^{(2)} \geq \cdots \geq 0).
\end{gather*}
Then we write $\lambda \trianglelefteq \mu$ if
\[ \sum_{i=1}^j \lambda_1^{(j)} \leq \sum_{i=1}^j \mu_1^{(j)} 
\qquad \mbox{for all $j$}\]
and 
\[ |\lambda_1|+ \sum_{i=1}^j \lambda_2^{(j)} \leq |\mu_1|+
\sum_{i=1}^j \mu_2^{(j)} \qquad \mbox{for all $j$}.\]
Note that, if $|\lambda_1|=|\mu_1|$, then we have 
\[ \lambda \trianglelefteq \mu \quad \Leftrightarrow \quad 
\lambda_1 \trianglelefteq \mu_1 \quad \mbox{and}\quad \lambda_2 
\trianglelefteq \mu_2,\]
where, on the right, the symbol $\trianglelefteq$ just denotes the
usual dominance order on partitions, as in Example~\ref{symgroup}.
The following result is a refinement  of \cite[Remark~3.7]{BI2}
(which only deals with elements of the same $t$-length).

\begin{prop} \label{twoprop1} Let $x,y\in W_n$ be such that $x\leq_{\cLR} y$. 
Assume that $x$ belongs to an RS-cell of shape $(\lambda_1,\lambda_2)$ 
and $y$ belongs to an RS-cell of shape $(\mu_1,\mu_2)$. Then we have 
\[ (\lambda_1,\lambda_2^*) \trianglelefteq (\mu_1,\mu_2^*),\]
with equality only if $x \sim_{\cLR} y$.
\end{prop}

\begin{proof} Let $w \in W_n$ and write $w=a_wa_l\sigma_wb_w^{-1}$
where $l=\ell_t(w)$. Now the parabolic subgroup $\fS_{l,n-l}$ is a direct 
product of $\fS_l=\langle s_1,\ldots,s_{l-1}\rangle$ and $\fS_{[l+1,n]}=
\langle s_{l+1},\ldots, s_{n-1} \rangle$. Thus, we have 
\[ \sigma_w=\sigma_w'\sigma_w'' \qquad \mbox{where $\sigma_w'\in 
\fS_l$ and $\sigma_w''\in \fS_{[l+1,n]}$}.\]
Since $\fS_{l}$ is a Coxeter group of type $A_{l-1}$, the classical 
Robinson--Schensted correspondence associates to $\sigma_w'$ a pair of 
tableaux whose shape is a partition of $l$, say $\nu_2$. Similarly, 
since $\fS_{[l+1,n]}$ is of type $A_{n-l}$, the classical 
Robinson--Schensted correspondence associates to $\sigma_w''$ a pair of 
tableaux whose shape is a partition of $n-l$, say $\nu_1$. Then, by
the discussion in \cite[4.7]{BI}, we have 
\begin{itemize}
\item[($\spadesuit$)] $w$ belongs to an RS-cell of shape $(\nu_1,
\nu_2^*)$.
\end{itemize}
Now consider the given two elements $x,y \in W_n$ such that $x\leq_{\cLR}
y$. By Proposition~\ref{mainbiI}, this implies $\ell_t(y)\leq \ell_t(x)$. In 
particular, $x$ and $y$ cannot lie in the same two-sided cell unless $x$ 
and $y$ have the same $t$-length. We now distinguish two cases.

\medskip
{\em Case 1.} Assume that $l:=\ell_t(x)=\ell_t(y)$. By Theorem~\ref{mainbi2}, 
we know that $x\leq_{\cLR} y$ implies that $\sigma_x \leq_{\cLR,l} \sigma_y$. 
Furthermore, it is well-known and easy to check that the Kazhdan--Lusztig 
pre-order relations are compatible with direct products; in particular, we 
have 
\[ \sigma_x' \leq_{\cLR}' \sigma_y' \qquad \mbox{and}\qquad 
\sigma_x'' \leq_{\cLR}'' \sigma_y'',\]
where a single dash denotes the pre-order relation on $\fS_l$ and a 
double-dash denotes the pre-order on $\fS_{[l+1,n]}$. Now ($\spadesuit$) 
shows that 
\begin{itemize}
\item $\sigma_x'$ is associated with the partition $\lambda_2^*$,
\item $\sigma_y'$ is associated with the partition $\mu_2^*$,
\item $\sigma_x''$ is associated with the partition $\lambda_1$,
\item $\sigma_y''$ is associated with the partition $\mu_1$.
\end{itemize}
Thus, we are reduced to statements concerning two-sided cells in the 
symmetric group. Now Example~\ref{symgroup}(e) shows that we have the 
implications
\[\sigma_x' \leq_{\cLR}' \sigma_y' \quad \Rightarrow \quad \lambda_2^* 
\trianglelefteq \mu_2^* \qquad \mbox{and}\qquad \sigma_x'' \leq_{\cLR}' 
\sigma_y'' \quad \Rightarrow \quad \lambda_1 \trianglelefteq \mu_1.\]
This yields $(\lambda_1, \lambda_2^*) \trianglelefteq (\mu_1,\mu_2^*)$
as required. Furthermore, if $(\lambda_1, \lambda_2^*)= (\mu_1, \mu_2^*)$, 
then $x\sim_{\cLR} y$ by Theorem~\ref{mainbi2}.

\medskip
{\em Case 2.} Assume that $\ell_t(y)<\ell_t(x)$. By Case~1 and the definition 
of $\leq_{\cLR}$, it is enough to consider an elementary step, where 
$\ell_t(y)< \ell_t(x)$ and $x\leftarrow_{\cL} y$ or $x \leftarrow_{\cR} y$.  
Since $w \sim_{\cLR} w^{-1}$ for all $w\in W_n$, we can even assume that 
$x\leftarrow_{\cR} y$, that is, $h_{y,s,x}\neq 0$ for some $s\in \{t,s_1, 
s_2,\ldots,s_{n-1}\}$. Since we are assuming $\ell_t(y)< \ell_t(x)$, we must 
have $s=t$ and $x=yt>y$ by Proposition~\ref{mainbiI}. Thus, it remains to 
consider the effect of multiplying with $t$ on the generalized 
Robinson--Schensted correspondence. We claim that, if $x=yt>y$, then
\begin{itemize}
\item $\lambda_1$ is obtained from $\mu_1$ by decreasing one part by $1$, and
\item $\lambda_2^*$ is obtained from $\mu_2^*$ by increasing one part by $1$.
\end{itemize}
This is seen as follows. Recall the basic ingredients of the generalized
Robinson-Schensted correspondence. We write $y\in W_n$ as a signed permutation
\[ \left(\begin{array}{cccc} 1 & 2 & \cdots & n \\ \varepsilon_1 \cdot p_1
& \varepsilon_2\cdot p_2 & \cdots & \varepsilon_n \cdot p_n\end{array}
\right) \qquad (\varepsilon_1=1)\]
where the sequence $p_1,\ldots,p_n$ is a permutation of $1,\ldots,n$ and
where $\varepsilon_i=\pm 1$ for all~$i$. The fact that $\varepsilon_1=1$
follows from our assumption that $yt>y$. Let $1\leq i_1 <i_2 <\cdots <i_k
\leq n$ be the sequence of indices where the sign is ``$+$'' and let $1\leq
j_1<j_2 <\cdots< j_l\leq n$ be the sequence of indices where the sign is
``$-$''. Applying the usual ``insertion algorithm'' to these two sequences
of numbers, we obtain a standard bitableau $(A^+(y),A^-(y))$ of size~$n$
and shape $(\mu_1,\mu_2)$, where $\ell_t(y)=|\mu_2|=l$. 

Now let us multiply $y$ on the right by $t$. Then the corresponding signed 
permutation is given by 
\[ \left(\begin{array}{cccc} 1 & 2 & \cdots & n \\ -\varepsilon_1 \cdot 
p_1 & \varepsilon_2\cdot p_2 & \cdots & \varepsilon_n \cdot p_n\end{array}
\right).\]
Thus, the only effect of multiplying by $t$ is to change the sign in the 
first position of the above array. Hence, in order to obtain the tableaux 
$A^+(yt)$ and $A^-(yt)$, we must apply the insertion algorithm to the 
sequences $p_{i_2}, \ldots,p_{i_k}$ and $p_{i_1},p_{j_1},p_{j_2},
\ldots,p_{j_l}$, respectively. Thus, we are reduced to a purely 
combinatorial statement. Using \cite[\S 5, Prop.~1]{Fulton}, one shows
that the partition giving the shape of $A^+(yt)$ is obtained from the 
partition giving the shape of $A^+(y)$ by decreasing one part by~$1$. 
The same argument shows that the partition giving the shape of $A^-(yt)$ 
is obtained from the partition giving the shape of $A^-(y)$ by increasing 
one part by~$1$. (A much more general statement can be found in  
\cite[Theorem~4.2]{Bhama}.)

Now the definition of $\trianglelefteq$ immediately shows that, 
if $(\lambda_1, \lambda_2)$ 
is obtained from $(\mu_1,\mu_2)$ by the above procedure, 
then $(\lambda_1, \lambda_2^*) \trianglelefteq (\mu_1,\mu_2^*)$ 
as required.
\end{proof}

\begin{cor} \label{twocor1} Let $x,y\in W_n$ be such that $x\leq_{\cLR} y$.
Then we have $\alpha_y \leq \alpha_x$, with equality only if $x\sim_{\cLR} y$.
\end{cor}

\begin{proof} Assume that $x$ belongs to an RS-cell of shape $(\lambda_1, 
\lambda_2)$ and $y$ belongs to an RS-cell of shape $(\mu_1,\mu_2)$. Hence, 
by Example~\ref{bn1}, we have 
\[ \alpha_x=\alpha_{(\lambda_1,\lambda_2^*)} \qquad \mbox{and}\qquad
\alpha_y=\alpha_{(\mu_1,\mu_2^*)}.\]
The description of the generalized Robinson--Schensted correspondence in
\cite{BI} shows that $\ell_t(x)=|\lambda_2|=|\lambda_2^*|$ and $\ell_t(y)=
|\mu_2|= |\mu_2^*|$. Now, by Proposition~\ref{mainbiI}, the condition 
$x\leq_{\cLR} y$ implies that $\ell_t(y)\leq \ell_t(x)$. 

{\em Case 1.} Assume that $l:=\ell_t(x)=\ell_t(y)$. Then we have
\[\lambda_1 \trianglelefteq \mu_1 \qquad \mbox{and}\qquad 
\lambda_2^* \trianglelefteq \mu_2^*\]
by Proposition~\ref{twoprop1}. Now note the following property of the 
dominance order. For any partitions $\nu$ and $\nu'$ of $n$, we have  
\[\nu^* \trianglelefteq \nu'^* \qquad \Leftrightarrow \qquad \nu' 
\trianglelefteq \nu \qquad \Rightarrow \qquad n(\nu) \leq n(\nu'),\]
with equality only for $\nu=\nu'$; see, for example, 
\cite[Exercise~5.6]{ourbuch}. Using the above property, we conclude that
\[ n(\mu_1)+2n(\mu_2^*)-n(\mu_2)\leq n(\lambda_1)+2n(\lambda_2^*)-
n(\lambda_2),\]
with equality only if $(\mu_1,\mu_2)=(\lambda_1,\lambda_2)$. Hence, the 
formula for $\alpha_w$ shows that $\alpha_y \leq \alpha_x$, as required. 
Furthermore, if $\alpha_x=\alpha_y$, then we necessarily have $(\lambda_1,
\lambda_2)=(\mu_1,\mu_2)$, and so $x\sim_{\cLR} y$ by Theorem~\ref{mainbi2}.

{\em Case 2.} Assume that $\ell_t(y)< \ell_t(x)$. As in the proof of
Proposition~\ref{twoprop1}, we can reduce to the case where $x=yt>y$.\\
In this case we have $\alpha_x-\alpha_y=b+a(2r'-m-r)$, where $m,r,r'$
are integers determined by the conditions : \[\lambda_1^{(m)}=\mu_1^{(m)}
-1,\quad \lambda_2^{(r)}=\mu_2^{(r)}+1,\quad \lambda_2^{*(r')}=
\mu_2^{*(r')}+1\, .\]
Now note that we have $2r'-m-r+n-1\geq 0$, hence $\alpha_x>\alpha_y$
as desired.
\end{proof}

\medskip
\noindent {\bf Proofs of Theorem~\ref{main}, Theorem~\ref{mainb} and
Theorem~\ref{Jring}.} Let us recall the principal ingredients. By 
Example~\ref{bn}, we know that $\bH_n$ is integral and normalized. 
Furthermore, by the discussion in Example~\ref{bn1}, the characters 
afforded by all left cells are irreducible. Finally, the assumption 
($\clubsuit$) in Lemma~\ref{cor1a} holds by Corollary~\ref{twocor1}. 
Thus, we can conclude that 
\[\ba(w)=\alpha_{(\lambda_1,\lambda_2)} \qquad \mbox{($w\in \cC_{(\lambda_1,
\lambda_2)}$)}.\]
The identification of $\cC_{(\lambda_1,\lambda_2)}$ in Example~\ref{bn1} 
and the formula for $\alpha_{(\lambda_1,\lambda_2)}$ in Example~\ref{bn} 
now yield the explicit description of the $\ba$-function of $W_n$, 
proving Theorem~\ref{main}. 

Now, once the $\ba$-function is determined, condition ($\clubsuit$) in
Lemma~\ref{cor1a} yields the implication ``$x \leq_{\cLR} y \Rightarrow 
\ba(y)\leq \ba(x)$'' for any $x,y \in W_n$, that is, {(P4)} holds. 
But then Corollary~\ref{twocor1} also yields the fact that, if 
$x\leq_{\cLR} y$ and $\ba(x)=\ba(y)$, then $x\sim_{\cLR} y$, that is,
{(P11)} holds. Now Lemma~\ref{cor1c} and Lemma~\ref{cor1d} show that all 
the other properties mentioned in Theorem~\ref{mainb} hold. As far as {(P12)}
is concerned, note that every parabolic subgroup of $W_n$ is a direct
product of a group of type $B_k$ and possibly some factors of type $A_{n_i}$.
Since {(P1)--(P15)} are known to hold for groups of type $A_{n_i}$, we
conclude that {(P3)}, {(P4)}, {(P8)} hold for every parabolic 
subgroup of $W_n$. This formally implies that {(P12)} holds; see 
\cite[14.12]{Lusztig03}.

Finally, as explained in \cite[18.3]{Lusztig03}, the ring $J$ can be
constructed once it is known that {(P1)--(P8)} are known to hold. The 
structure of $J$ is now determined by Proposition~\ref{strucJ} (and
Remark~\ref{strucJ1}). 

\begin{rem} \label{bonnaf} Bonnaf\'e has remarked at the end of
\cite[\S 4]{BI2} that, once the equality $\ba(z)=\alpha_z$ ($z\in W_n$) 
is established, the methods in his paper \cite{BI2} yield properties 
(P1), (P4), (P6), (P11), (P12) and the first assertion of (P13). However,
it does not seem to be possible to gain control over the constants
$\gamma_{x,y,z}$ in his approach, and this is where the leading matrix 
coefficients of orthogonal representations naturally come into play.
\end{rem}

\begin{exmp} \label{expb2} Let us consider the case $n=2$, where
$W_2=\langle t,s_1\rangle$ is the dihedral group of order~$8$. 
We set $s_0=t$. The polynomials $P_{y,w}^*$ and the left cells have 
already been determined by an explicit computation in 
\cite[\S 6]{Lusztig83}. The left cells are
\[ \{1\},\quad \{s_1\}, \quad \{s_0,s_1s_0\}, \quad \{s_0s_1,s_1s_0s_1\},
\quad \{s_0s_1s_0\}, \quad \{w_0\}\]
where $w_0=s_1s_0s_1s_0$ is the unique element of maximal length.
Using the polynomials $P_{y,w}^*$, we compute:
\begin{alignat*}{2}
\bD(1) &= 0, &\qquad n_1&=1,\\
\bD(s_1) &= a, &\qquad n_{s_1}&=1,\\
\bD(s_0) &= b, &\qquad n_{s_0}&=1,\\
\bD(s_1s_0s_1) &= b, &\qquad n_{s_1s_0s_2}&=1,\\
\bD(s_0s_1s_0) &= 2b-a, &\qquad n_{s_0s_1s_0}&=-1,\\
\bD(w_0) &= 2b+2a, &\qquad n_{w_0}&=1.
\end{alignat*}
In particular, we see that the coefficients $n_z$ can be negative.
Proposition~\ref{strucJ} shows that 
\[ t_{s_0s_1s_0}^2=-t_{s_0s_1s_0}; 
\qquad \mbox{see also Lusztig \cite[18.7]{Lusztig03}}.\]
There is a unique irreducible character $\chi_\lambda$ of degree~$2$; it 
is labelled by $\lambda=((1),(1))$. A corresponding orthogonal 
representation and the leading matrix coefficients are explicitly 
described in \cite[Exp.~5.5]{my02}. We have $\alpha_\lambda=b$ and the 
representation $\fX_\lambda$ is given by 
\[\fX_\lambda \colon T_{s_0} \mapsto\begin{bmatrix} V & 0\\0 & -V^{-1}
\end{bmatrix}, \quad T_{s_1} \mapsto \frac{1}{V^2+1}\begin{bmatrix}
v-v^{-1} & 1+V^2v^{-2}\\ 1+V^2v^2 & V^2(v-v^{-1})\end{bmatrix}.\]
Thus, we see that $c_{s_0,\lambda}^{2,2}=c_{s_1s_0s_1,\lambda}^{1,1} =1$
and $c_{s_0s_1,\lambda}^{2,1}=c_{s_1s_0,\lambda}^{1,2}=-1$.
We conclude that 
\[ s_1s_0s_1=w_\lambda(1,1),\quad s_1s_0=w_\lambda(1,2),\quad
s_0s_1=w_\lambda(2,1),\quad s_0=w_\lambda(2,2).\]
Hence, Proposition~\ref{strucJ} yields  a ring isomorphism $J_\lambda 
\cong M_2(Z)$ where 
\[t_{s_1s_0s_1} \mapsto \begin{bmatrix} 1 & 0 \\ 0 & 0 \end{bmatrix},\quad
t_{s_1s_0} \mapsto \begin{bmatrix} 0 & -1 \\ 0 & 0 \end{bmatrix},\quad
t_{s_0s_1} \mapsto \begin{bmatrix} 0 & 0 \\ -1 & 0 \end{bmatrix},\quad
t_{s_0} \mapsto \begin{bmatrix} 0 & 0 \\ 0 & 1 \end{bmatrix}.\]
Note that one also obtains an isomorphism $J_\lambda \cong M_2(\Z)$ by 
sending the above four elements $t_x$ directly to the corresponding matrix
units (omitting the sign in the matrices), as in \cite[18.7]{Lusztig03}. The
signs arise from the construction in the proof of Proposition~\ref{strucJ} 
and the choice of the orthogonal representation. Note that the latter is 
not unique: for example, one can conjugate $\fX_\lambda$ by a diagonal 
matrix with $\pm 1$ on the diagonal.
\end{exmp}

\medskip
\noindent {\bf Proof of Theorem~\ref{murphy}.} Let us recall some 
ingredients of the construction of the Dipper--James--Murphy basis. Let 
$\lambda=(\lambda_1, \lambda_2)$ be a bipartition of $n$; let $(\fs,\ft)$ 
be a pair of standard bitableaux of shape $\lambda$. By \cite[\S 4]{DJM}, 
we have $x_{\fs\ft}=T_dx_\lambda T_{d'}$ where $d,d'$ are certain elements 
in $\Sym_n$; the element $x_\lambda$ is defined in \cite[4.1]{DJM}. 
Since $N^\lambda$ is a two-sided ideal, we conclude that 
\[ N^\lambda=\sum_{\mu\in \Lambda_n;\, \lambda \trianglelefteq \mu}
\bH_n x_\mu \bH_n. \]
We must show that $N^\lambda=M^\lambda$, where we set
\[ M^\lambda:=\Big\langle C_w' \,\Big|\begin{array}{l} \mbox{ $w$ corresponds 
to an RS-cell of shape}\\ \mbox{ $\nu=(\nu_1,\nu_2)$ where $(\lambda_1,
\lambda_2) \trianglelefteq (\nu_2,\nu_1^*)$}\end{array}\Big\rangle_A
\subseteq \bH_n.\]
One easily checks that $(\nu_1,\nu_2^*) \trianglelefteq (\lambda_2^*,
\lambda_1^*) \Leftrightarrow (\lambda_1,\lambda_2) \trianglelefteq 
(\nu_2,\nu_1^*)$ (using the analogous statement for the dominance order 
on partitions; see \cite[Exercise~5.6]{ourbuch}). Hence, by 
Proposition~\ref{twoprop1}, $M^\lambda$ is a two-sided ideal of $\bH_n$. 

We now show that $N^\lambda\subseteq M^\lambda$.
Let $\mu=(\mu_1,\mu_2)$ be a bipartition of $n$ such that $(\lambda_1,
\lambda_2) \trianglelefteq (\mu_1,\mu_2)$. Let $l:=|\mu_1|$. The element 
$x_\mu$ is defined as the product of three factors $u_l^+$, $x_{\mu_1}$, 
$x_{\mu_2}$. The formula in Bonnaf\'e \cite[Prop.~2.5]{BI2} shows that, 
up to multiplying by a monomial in $V$ and $v$, the factor $u_l^+$ equals 
$T_{\sigma_l}\, C_{a_l}'$ where $\sigma_l$ is the longest element in 
$\fS_l=\langle s_1,\ldots,s_{l-1} \rangle$. Furthermore, by Lusztig 
\cite[Cor.~12.2]{Lusztig03}, we have $x_{\mu_1} x_{\mu_2}=
v^{l(w_\mu)}C_{w_\mu}'$ where $w_\mu$ is the longest element in the 
Young subgroup of $\fS_{l,n-l}$ given by $\mu=(\mu_1,\mu_2)$. Finally, by 
\cite[Prop.~2.3]{BI2}, we have $C_{a_l}' C_{w_\mu}'=C_{a_l w_\mu}'$. 
Hence, we obtain 
\[ x_\mu=T_{\sigma_l}\,  C_{a_lw_\mu}'\qquad \mbox{(up to multiplying 
by a monomial in $V$ and $v$)}.\]
By relation ($\spadesuit$) in the proof of Proposition~\ref{twoprop1}, 
$a_lw_\mu$ belongs to an RS-cell of shape $(\mu_2^*,\mu_1)$. Hence, 
since $(\lambda_1, \lambda_2) \trianglelefteq (\mu_1,\mu_2)$ and since 
$M^\lambda$ is an ideal, we obtain that $x_\mu\in M^\lambda$. As this 
holds for all $\mu$ such that $\lambda \trianglelefteq \mu$, we conclude 
that $N^\lambda \subseteq M^\lambda$. 

In order to show equality, we note that $M^\lambda$ is free over 
$A$ of rank $\sum_{\lambda \trianglelefteq \nu} d_\nu^2$, where $d_\nu$ 
denotes the number of standard bitableaux of shape~$\nu$. By \cite[4.15]{DJM}, 
$N^\lambda$ is free over $A$ of the same rank. Consequently, we have 
$K_0\otimes_A N^\lambda =K_0 \otimes_A M^\lambda$, where $K_0$ is the field
of fractions of $A$. So there exists some $0\neq f \in A$ such that 
$fM^\lambda \subseteq N^\lambda \subseteq M^\lambda$. Now, since the 
generators of $N^\lambda$ can be extended to an $A$-basis of $\bH_n$
(see \cite[\S 4]{DJM}), the quotient $\bH_n/N^\lambda$ is a free $A$-module. 
Hence, $M^\lambda \subseteq N^\lambda$ and the conclusion follows. 

\medskip
\noindent {\bf Acknowledgements.} This paper was written while the first
named author enjoyed the hospitality of the Bernoulli Center at the EPFL 
Lausanne (Switzerland), during the research program ``Group representation 
theory'' from january to june 2005. The second named author gratefully 
acknowledges support by the Fonds National Suisse de la Recherche
Scientifique.


\end{document}